\newif\ifArxivFormat

\ArxivFormattrue

\ifArxivFormat
	\documentclass[onefignum,onetabnum,a4paper,oneside]{siamart190516}
	\usepackage[top=1in, bottom=0.75in, left=0.75in, right=0.75in]{geometry}
\else
	\documentclass[review,onefignum,onetabnum]{siamart190516}
\fi

\usepackage{lipsum}
\usepackage{amsfonts}
\usepackage{graphicx}
\usepackage{epstopdf}
\usepackage{amsopn}
\usepackage{mymath}
\usepackage{arydshln} % dashed line in array
\usepackage{parskip}
\usepackage{nth}
\usepackage{varwidth, tikz}
\usetikzlibrary{shapes, arrows, shapes.misc, arrows.meta, positioning, matrix, calc, fit, fadings,patterns}
\usepackage{csvsimple}
\begingroup
\definecolorseries{gTob}{rgb}{last}{black!10!white}{black!90!white}
\endgroup
\usepackage{pgfplots}
\pgfplotsset{compat=newest}
\usepgflibrary{plothandlers,plotmarks}
\newbox\aMark
\setbox\aMark\hbox{\begin{pgfpicture}\pgfuseplotmark{*}\end{pgfpicture}}
\newbox\bMark
\setbox\bMark\hbox{\begin{pgfpicture}\pgfuseplotmark{triangle*}\end{pgfpicture}}
\newbox\cMark
\setbox\cMark\hbox{\begin{pgfpicture}\pgfuseplotmark{square*}\end{pgfpicture}}

\usepackage{graphicx,epstopdf} % <- Preamble
\usepackage{subcaption}
\usepackage{graphicx,xcolor} 
\usepackage[framemethod=tikz]{mdframed}

\usetikzlibrary{external}
\ifpdf
  \DeclareGraphicsExtensions{.eps,.pdf,.png,.jpg}
\else
  \DeclareGraphicsExtensions{.eps}
\fi
\usepackage{booktabs}
\usepackage{siunitx}

\usepackage{algorithm}
\usepackage{algorithmic}

\newcommand{\cplxity}[1]{\textcolor{gray!50}{$\mathcal{O}\left(#1\right)$}}
\usepackage{letltxmacro}
\newlength{\commentindent}
\makeatletter
\renewcommand{\algorithmiccomment}[1]{\unskip\hfill\makebox[\commentindent][l]{\footnotesize{\cplxity{#1}}}\par}
\LetLtxMacro{\oldalgorithmic}{\algorithmic}
\renewcommand{\algorithmic}[1][0]{%
  \oldalgorithmic[#1]%
  \renewcommand{\ALC@com}[1]{%
    \ifnum\pdfstrcmp{##1}{default}=0\else\algorithmiccomment{##1}\fi}%
}
\makeatother

\newsiamremark{remark}{Remark}
\newsiamremark{hypothesis}{Hypothesis}
\crefname{hypothesis}{Hypothesis}{Hypotheses}
\newsiamthm{claim}{Claim}

\headers{A Higher-Order GSVD for Rank Deficient Matrices}{I. Kempf, P. J. Goulart and S. R. Duncan}

\ifArxivFormat
\title{A Higher-Order Generalized Singular Value Decomposition for Rank Deficient Matrices\thanks{Submitted \today.
\funding{The research leading to these results was supported in part by the Diamond Light Source and in part by the Engineering and Physical Sciences Research Council (EPSRC) under a Collaborative Awards in Science and Engineering (CASE) studentship.}}}
\else
\title{A Higher-Order Generalized Singular Value Decomposition for Rank Deficient Matrices\thanks{Submitted to the editors DATE.
\funding{The research leading to these results was supported in part by the Diamond Light Source and in part by the Engineering and Physical Sciences Research Council (EPSRC) under a Collaborative Awards in Science and Engineering (CASE) studentship.}}}
\fi

\author{Idris Kempf\thanks{Department of Engineering Science, University of Oxford, Parks Road, Oxford, OX1 3PJ, UK
  (\email{idris.kempf@eng.ox.ac.uk},\email{paul.goulart@eng.ox.ac.uk}, \email{stephen.duncan@eng.ox.ac.uk}).}
\and Paul J. Goulart\footnotemark[2]
\and Stephen R. Duncan\footnotemark[2]}

\newcommand{\phantomproofbox}{{\color{white}\vbox{\hrule height0.6pt\hbox{\vrule height1.3ex width0.6pt\hskip0.8ex\vrule width0.6pt}\hrule height0.6pt}}}
\theoremstyle{nonumberplain}
\theoremheaderfont{\normalfont\itshape}
\theorembodyfont{\normalfont}
\theoremseparator{.}
\theoremsymbol{\phantomproofbox}

\usepackage{xcolor,calc}

\ifpdf
\hypersetup{
  pdftitle={A Higher-Order Generalized Singular Value Decomposition for Rank-Deficient Matrices},
  pdfauthor={I. Kempf, P. Goulart and S. Duncan}
}
\fi

\usepackage{xcolor}
\definecolor{colorblue}{RGB}{0,33,71}  		% blue
\ifArxivFormat
\newcommand{\added}[1]{#1}
\else
\newcommand{\added}[1]{\textcolor{red!90}{#1}}
\fi

\usepackage{enumerate}
\usepackage[shortlabels]{enumitem}

\newsiamremark{newremark}{\added{Remark}}
\newsiamthm{newlemma}{\added{Lemma}}
\newsiamthm{newdefinition}{\added{Definition}}
\newsiamthm{newthm}{\added{Theorem}}
\newsiamthm{newcorollary}{\added{Corollary}}

\begin{document}

\maketitle

% REQUIRED
\begin{abstract}
The higher-order generalized singular value decomposition (HO-GSVD) is a matrix factorization technique that extends the GSVD to $N \ge 2$ data matrices, and can be used to identify common subspaces that are shared across multiple large-scale datasets with different row dimensions. The standard HO-GSVD factors $N$ matrices $\inR{A_i}{m_i}{n}$ as $A_i=U_i\Sigma_i\trans{V}$, but requires that each of the matrices $A_i$ has full column rank. We propose a modification of the HO-GSVD that extends its applicability to rank-deficient data matrices $A_i$. If the matrix of stacked $A_i$ has full rank, we show that the properties of the original HO-GSVD extend to our approach. We extend the notion of common subspaces to isolated subspaces, which identify features that are unique to one $A_i$. We also extend our results to the higher-order cosine-sine decomposition (HO-CSD), which is closely related to the HO-GSVD. Our extension of the standard HO-GSVD allows its application to datasets with with $m_i<n$ or $\rank(A_i)<n$, such as are encountered in bioinformatics, neuroscience, control theory or classification problems.
\end{abstract}

% REQUIRED

% REQUIRED

\section{Introduction}\label{sec:introduction}
The \emph{generalized singular value decomposition} (GSVD) \cite{LOAN_GSVD} is an extension of the well-known \emph{singular value decomposition} (SVD) to $N=2$ matrices. The GSVD decomposes a pair of matrices $\inR{A_1}{m_1}{n}$ and $\inR{A_2}{m_2}{n}$, with $m_1\geq n$, by factorizing each of the matrices as $A_i=U_i\Sigma_i\trans{V}$. The matrix of right generalized singular vectors $\inR{V}{n}{n}$, with $\det(V)\neq 0$, is shared between the decompositions, but unlike the standard SVD is not an orthogonal matrix. The columns of the matrices $\inR{U_i}{m_i}{m_i}$ are commonly referred to as \emph{left generalized singular vectors}, each satisfying $\xTx{U_i}=I$. The matrices $\Sigma_i=\diag(\sigma_{i,1},\dots,\sigma_{i,r})$, with $\sigma_{i,k}\geq 0$  and $r=\rank([\trans{A}_1,\,\trans{A}_2\,\trans{]})$, contain the \emph{generalized singular values}. \added{The generalized singular values measure the significance of the right generalized singular vectors $v_k$ in the factorization of each $A_i$~\cite{HOGSVD}. If $\sigma_{1,k}=\sigma_{2,k}$, then $v_k$ solves the \emph{generalized singular value problem} $\xTx{A_1}v_k=\mu\xTx{A_2} v_k$ with $\mu=1$~\cite[Ch. 8.7]{GOLUB4}.}

The \emph{higher-order GSVD} (HO-GSVD)~\cite{HOGSVD} is an extension of the GSVD to $N\geq 2$ matrices. Given $N$ matrices $A_1,\dots,A_N$, the HO-GSVD decomposes each $A_i$ as
\begin{align}\label{eq:factorization}
A_i=U_i\Sigma_i\trans{V}, \qquad i = 1,\dots,N,
\end{align}
where $\inR{U_i}{m_i}{n}$, $\inR{\Sigma_i}{n}{n}$ and $\inR{V}{n}{n}$ with $\det(V)\neq 0$ being shared among all factorizations. The matrix $V$ is obtained from the eigensystem ${S_\pi}V=V\varSigma$, where $\varSigma\eqdef\diag(\varsigma_1,\dots,\varsigma_n)$ and $S_\pi$ is the \emph{arithmetic mean of all pairwise quotients} $D_{i,\pi}\inv{D_{j,\pi}}$,
\begin{align}\label{eq:S}
S_\pi\eqdef\frac{1}{N(N-1)}\sum_{i=1}^N\sum_{j=i+1}^N\left(D_{i,\pi}\inv{D_{j,\pi}}+D_{j,\pi}\inv{D_{i,\pi}}\right),
\end{align}
with $D_{i,\pi}$ defined as
\begin{align}\label{eq:Di}
D_{i,\pi}\eqdef\xTx{A_i}+\pi \xTx{A},
\qquad\pi\geq 0,
\end{align}
where $A\eqdef [\trans{A}_1,\dots,\trans{A}_N\trans{]}$. \added{The case $\pi=0$ corresponds to the standard HO-GSVD framework~\cite{HOGSVD}. Under the assumption that $A$ has full column rank,} introducing the term $\pi \xTx{A}$ with $\pi>0$ will allow us to accommodate rank-deficient matrices $A_i$.

\added{%%%
Using the factorization~\eqref{eq:factorization}, the matrices $A_i$ can be rewritten as
\begin{align}\label{eq:factorization2}
A_i = \underbrace{\sum_{k\in\mathcal{I}_N}\sigma_{i,k}u_{i,k}\trans{v_k}}_{\text{common}} +
	  \underbrace{\sum_{k\in\mathcal{I}_1}\sigma_{i,k}u_{i,k}\trans{v_k}}_{\text{isolated}} +
	  \sum_{k\in\mathcal{I}_\perp}\sigma_{i,k}u_{i,k}\trans{v_k},
\end{align}
where $\mathcal{I}_N\,\cup\,\mathcal{I}_1\,\cup\,\mathcal{I}_\perp=\lbrace 1,\dots,n\rbrace$ and $\mathcal{I}_N$, $\mathcal{I}_1$ and $\mathcal{I}_\perp$ are mutually disjoint.} The columns $u_{i,k}$ of the matrices $U_i$ are referred to as \emph{left basis vectors}, and the diagonal matrices $\Sigma_i$ contain the generalized singular values $\sigma_{i,k}$. The \emph{right basis vectors} $v_k$ are shared across all $A_i$. It can be shown that the GSVD is a special case of the HO-GSVD with $N=2$ and that the standard SVD of $A_j$ can be obtained from setting $A_i=\I$ for $i\neq j$ and~$N\geq 2$~\cite{HOGSVD}. For the case $\pi = 0$, it was shown in \cite{HOGSVD} that the subspace associated with the unit eigenvalues of $S_\pi$ forms the \emph{common HO-GSVD subspace} (see Section~\ref{sec:commonsubspaces} and Definition~\ref{def:commonsubspace}), which is preserved for $\pi > 0$. \added{This subspace is spanned by the right basis vectors ${v_k}$, $k\in\mathcal{I}_N$, for which $\sigma_{i,k}=\sigma_{j,k}$, and the associated left generalized singular vectors $u_{i,k}$ are orthogonal to $u_{i,j}$, $j\neq k$.}

The HO-GSVD is a technique that is of particular use in multimodal data fusion~\cite{MMDFUSION}, which aims to identify common features across multiple data sets that describe related phenomena. \added{Many tensor or multi-matrix decompositions are obtained from extending single-matrix factorizations to multiple matrices, such as the parallel factor analysis (PARAFAC~\cite{PARAFAC} or PARAFAC2~\cite{PARAFAC2}), multilinear SVDs~\cite{HOSVD}, multilinear principal component analysis~\cite{MPCA} or the higher-order eigenvalue decomposition~\cite{HOEVD}. The different extensions preserve some but not all of the single-matrix factorization properties~\cite{HOGSVDTHESIS}, such as exactness, orthogonality or rank conditions of the factor matrices.} Some tensor decompositions require that the matrices $A_i$ share the same dimensions, e.g. a third-order tensor $\mathcal{A}=A_1\times A_2\times\dots A_N$ requires that all matrices $A_i$ have dimensions $m\times n$, which imposes constraints on the data acquisition. In contrast, the HO-GSVD is an exact matrix factorization so that $A_i=U_i\Sigma_i\trans{V}$ for $i=1,\dots,N$, and it can accommodate $\inR{A_i}{m_i}{n}$ with different row dimensions $m_i$, although no constraints, such as orthogonality, can be imposed on the factor matrices.

One shortcoming of the original HO-GSVD framework~\cite{HOGSVD} is that the arithmetic mean~\eqref{eq:S} is only well defined for matrices $A_i$ that have full column rank. If $\rank(A_i)<n$ or $m_i<n$ for some $i$, then the inverse $\inv{(\xTx{A_i})}$ does not exist and so $S_0$ in~\eqref{eq:S} is not well defined. In addition, computing~\eqref{eq:S} may be inaccurate when one or more of the $A_i$ have small singular values. 

Provided that the matrix $A$ of stacked $A_i$ has full rank, introducing the term $\pi\xTx{A}$ in~\eqref{eq:Di} with parameter $\pi>0$ has the effect of shifting the eigenvalues of each $D_{i,\pi}$, so that the terms $D_{i,\pi}$ are guaranteed to be invertible and the HO-GSVD can be computed for $A_i$ with arbitrary rank. When all $A_i$ have full column rank, we show that $S_\pi$ with $\pi>0$ and $S_0$ both capture the common subspaces of $A_1,\dots,A_N$. We introduce the notion of an \emph{isolated HO-GSVD subspace} that accounts for the fact that a rank deficient $A_i$ can have a non-empty (right) nullspace. The isolated HO-GSVD subspace is spanned by the right basis vectors ${v_k}$, $k\in\mathcal{I}_1$, for which $\sigma_{i,k}> 0$ and $\sigma_{j,k}=0$, $j\neq i$. The associated left basis vectors $u_{i,k}$ are orthogonal to $u_{i,l}$, $l\neq k$, $i=1,\dots,N$.

The GSVD is closely related to the (thin) \emph{cosine-sine decomposition} (CSD)~\cite[Ch. 2.5.4]{GOLUB4}. In essence, the CSD states that the SVDs of $\inR{Q_1}{m_1}{n}$ and $\inR{Q_2}{m_2}{n}$ satisfying $\xTx{Q_1}+\xTx{Q_2}=\I$ share the same matrix of standard right singular vectors~\cite{LOAN_COMP_GSVD}. The GSVD can be obtained from applying a CSD to the matrices $Q_1$ and $Q_2$ that are obtained from the thin QR factorization of the stacked matrices $[\trans{A}_1,\,\trans{A}_2\trans{]}=QR$, where $Q$ is conformably partitioned such that $A_i=Q_iR$.

Analogous to the GSVD and the CSD, the HO-GSVD is closely related to the \emph{higher-order CSD} (HO-CSD)~\cite{LOANSLIDES}. The HO-GSVD of $N$ matrices $A_i$ can be obtained from the HO-CSD of $Q_1,\dots,Q_N$ that are obtained from the thin QR factorization of the stacked matrices $[\trans{A}_1,\,\dots,\,\trans{A}_N\trans{]}$. As in the case of the HO-GSVD, the computation of the HO-CSD proposed in~\cite{LOANSLIDES} is limited to the case that all $Q_i$ have full rank. In this paper, we also propose to compute the HO-CSD in a different way, which allows for the factorization of rank-deficient $Q_i$ satisfying $\xTx{Q_1}+\dots+\xTx{Q_N}=\I$.

The paper is organized as follows. Section~\ref{sec:main} presents the HO-CSD and the HO-GSVD that are applicable to rank-deficient matrices. In Section~\ref{sec:commonsubspaces}, we extend the notion of common HO-CSD and HO-GSVD subspaces to rank-deficient matrices. The effect of the parameter $\pi$ is investigated in Section~\ref{sec:pi}, followed by relating our findings to existing methods in Section~\ref{sec:comparison}. In Section~\ref{sec:computation}, we propose an algorithm for computing the HO-GSVD and the isolated subspace. The paper is concluded with an example application of the HO-GSVD in Section~\ref{sec:application}.

We use standard notation throughout the paper with $\range(A)$ and $\kernel(A)$ denoting the range and kernel of a matrix $A$. Positive-definite and positive-semidefinite matrices are denoted by $A\succ 0$ and $A\succeq 0$, respectively, and $\Rpos$ denotes the set of strictly positive real numbers.

\section{Main results}\label{sec:main}
Given $N$ matrices $\inR{A_i}{m_i}{n}$, let $A$ denote the matrix of stacked $A_i$ and $QR=A$ its thin QR factorization,
\begin{align}\label{eq:A}
A=
\begin{bmatrix}
A_1\\[-0.5em] \vdots\\[-0.1em] A_N
\end{bmatrix} = QR=
\begin{bmatrix}
Q_1\\[-0.5em] \vdots\\[-0.1em] Q_N
\end{bmatrix}
R,\qquad\inR{Q_{i}}{m_i}{n},\qquad\inR{R}{n}{n},
\end{align}
where it holds that
\added{\begin{align}\label{eq:assumptionQ}
&\xTx{Q}=\sum_{i=1}^N\xTx{Q_i}=\I,&\twonorm{Q_i}\leq 1\,\,\forall i=1,\dots,N.
\end{align}}
The matrices $A_i=Q_iR$ can individually have arbitrary rank, but \added{throughout the paper it is assumed that
\begin{align}\label{eq:assumptionA}
\rank(A)=\rank
\begin{pmatrix}
A_1\\[-0.5em] \vdots\\[-0.1em] A_N
\end{pmatrix}
=n,
\end{align}
so that $\det(R)\neq 0$ and $M\eqdef\sum_{i=1}^n m_i \geq n$. If~\eqref{eq:assumptionA} does not hold, the matrix $A$ can be padded using an additional matrix $A_{N+1}$ (see Remark~\ref{rem:rankdeficient}).} The quotient terms $D_{i,\pi}$~\eqref{eq:Di2} of the arithmetic mean $S_\pi$~\eqref{eq:S} can be rewritten as
\begin{align}\label{eq:Di2}
D_{i,\pi}=\xTx{A_i} + \pi\xTx{A} = \trans{R}\left(\xTx{Q_i} + \pi\I\right)R,
\end{align}
with parameter $\pi>0$. Since $\xTx{A_i}\succeq 0$ and $\pi\xTx{A}\succ 0$, the terms $D_{i,\pi}$ are guaranteed to be invertible.

The majority of our developments are based on the HO-CSD. Define $T_\pi$ as
\begin{align}\label{eq:T}
T_\pi \eqdef \frac{1}{N}\sum_{i=1}^N\invbr{\xTx{Q_i}+\pi\I},
\end{align}
\added{where it is assumed that~\eqref{eq:assumptionQ} holds.} The eigensystem of $T_\pi$ leads to the HO-CSD of the matrices $Q_i$. It can be shown (Appendix~\ref{app:SandT}) that $S_\pi$ and $T_\pi$ are related by:
\begin{align}\label{eq:SandT}
\tinv{R}S_\pi\trans{R} = \frac{1}{N-1}\left(\left(1+\pi N\right)T_\pi-\I\right).
\end{align}
\begin{theorem}\label{thm:T}
\added{Let $T_\pi$ be defined by~\eqref{eq:T} and suppose that~\eqref{eq:assumptionQ} holds.} There exists an orthogonal $\inR{Z}{n}{n}$ such that
\begin{equation}
\trans{Z}T_{\pi}Z=\diag(\tau_1,\dots,\tau_n), \label{thm:T:Zdef}
\end{equation}
where the columns of $Z$ are eigenvectors of $T_{\pi}$ and the eigenvalues $\tau_i$ of $T_{\pi}$ satisfy \[\tau_i\in\left[\tau_\text{min},\tau_\text{max}\right]\eqdef\left[ \invbr{\inv{N}+\pi},\,\frac{N-1}{N}\inv{\pi}+\frac{1}{N}\invbr{1+\pi}\right].
\]
\end{theorem}
\added{For the proof of Theorem~\ref{thm:T}, the following lemma is used.}
\ifArxivFormat
\begin{lemma}\label{thm:P}
\else
\begin{newlemma}\label{thm:P}
\fi
Let $\inR{P=\trans{P}}{n}{n}$ with $0\preceq P\preceq I$. For all $t\in\R^n$ with $\twonorm{t}=1$ and $\pi\geq 0$, it holds that $\trans{t}\left(\pi(1+\pi)\invbr{\pi I + P}\right)t\leq \trans{t}\left((1+\pi)I - P\right)t$. Moreover, equality holds iff $P$ has $p\geq 1$ eigenvalues $\lambda_1,\dots,\lambda_p\in\lbrace 0,1\rbrace$ associated with eigenvectors $v_1,\dots,v_p$, and $t\in\spanv{v_1,\dots,v_p}$.
\ifArxivFormat
\end{lemma}
\else
\end{newlemma}
\fi
\ifArxivFormat
\begin{proof}
\else
\begin{newproof}
\fi
The inequality $\trans{t}\left(\pi(1+\pi)\invbr{\pi I + P}\right)t\leq \trans{t}\left((1+\pi)I - P\right)t$ holds iff
\begin{align}\label{eq:ineqP}
(1+\pi)I - P - \pi(1+\pi)\invbr{\pi I + P} \succeq 0.
\end{align}
Set $P=V\Lambda\trans{V}$ with $\Lambda=\diag(\lambda_1,\dots,\lambda_n)$, $\lambda_i\in[0,1]$, so that~\eqref{eq:ineqP} amounts to \[f_i(\lambda_i)\eqdef 1+\pi-\lambda_i-\frac{\pi(1+\pi)}{\pi+\lambda_i}\geq 0,\qquad\qquad i=1,\dots,n.\]
Since $f_i''(\lambda_i)=-2\pi(1+\pi)/(\pi+\lambda_i)^3<0$ for $\lambda_i\in[0,1]$, the function $f_i(\lambda_i)$ is concave on $\lambda_i\in[0,1]$ and hence $f_i(\lambda_i)\geq\min\lbrace f_i(0),f_i(1)\rbrace=\min\lbrace 0,0\rbrace =0$ $\forall i=1,\dots,n$. Equality therefore holds iff $\lambda_i\in\lbrace 0,1\rbrace$.

For the second part of the claim, set $t=Va$ with $\twonorm{a}=1$, and pre- and post-multiply~\eqref{eq:ineqP} with $\trans{t}$ and $t$, respectively, to obtain
\begin{align}\label{eq:ineqfa}
\sum_{i=1}^n f_i(\lambda_i)a_i^2\geq 0.
\end{align}
Suppose that $t\in\spanv{v_1,\dots,v_p}$, then $\sum_{i=1}^n f_i(\lambda_i)a_i^2=\sum_{i=1}^p f_i(\lambda_i)a_i^2=0$. For the converse, suppose that $t\not\in\spanv{v_1,\dots,v_p}$ and that equality holds in~\eqref{eq:ineqfa}. Then there exists $j\in\lbrace p+1,\dots,n\rbrace$ with $a_j > 0$ and $f_j(\lambda_j) >0$, which is a contradiction.
\ifArxivFormat
\end{proof}
\else
\end{newproof}
\fi
\ifArxivFormat
\begin{proof}[\added{Proof of Theorem}~\ref{thm:T}]
\else
\begin{newproof}[\added{Proof of Theorem}~\ref{thm:T}]
\fi
The existence of a matrix $\inR{Z}{n}{n}$, $\xTx{Z}=I$, that diagonalizes $T_{\pi}$ is a consequence of the symmetry in \eqref{eq:T}. For the lower bound, substitute $u=(\trans{Q}_iQ_i+\pi I)^{\frac{1}{2}}t$ and $v=(\trans{Q}_iQ_i+\pi I)^{\sm\frac{1}{2}}t$ with $\twonorm{t}=1$ in the Cauchy-Schwarz inequality $(\trans{u}v)^2\leq\twonorm{u}^2\twonorm{v}^2$ to obtain
\begin{align}\label{eq:cauchyineq}
\trans{t}\invbr{\xTx{Q_i}+\pi I}t\geq \invbr{\trans{t}(\xTx{Q_i}+\pi I)t}.
\end{align}
Using~\eqref{eq:cauchyineq} and the harmonic-mean arithmetic-mean (HM-AM) inequality~\cite[Thm. 16]{MEANS}, a lower bound on $\trans{t}T_\pi t$ can be established as
\begin{subequations}\label{eq:ineqTT}\begin{align}
\trans{t}T_\pi t = \frac{1}{N}\sum_{i=1}^N\trans{t}\invbr{\xTx{Q_i}+\pi I}t
&\geq\frac{1}{N}\sum_{i=1}^N\frac{1}{\trans{t}(\xTx{Q_i}+\pi I)t}\label{eq:ineqTa}\\
&\geq\frac{N}{\pi N+\sum_{i=1}^N\trans{t}(\xTx{Q_i})t} = \tau_\text{min}\label{eq:ineqTb}.
\end{align}\end{subequations}
For the upper bound, apply Lemma~\ref{thm:P} with $P=\xTx{Q_i}$ to each summand of $T_\pi$:
\begin{align}\label{eq:ineqTupper}
\trans{t}T_\pi t\leq\frac{1}{N}\sum_{i=1}^N\trans{t}\left(\frac{1}{\pi}I-\frac{1}{\pi(1+\pi)}\xTx{Q_i}\right)t
= \frac{1}{\pi} - \frac{1}{N\pi(1+\pi)} = \tau_\text{max}.
\end{align}
\ifArxivFormat
\end{proof}
\else
\end{newproof}
\fi
%%%%%%%%%%%%%%%%%%%%%%%%

\begin{theorem}\label{thm:S}
\added{Let $S_\pi$ be defined by~\eqref{eq:S} and suppose that~\eqref{eq:assumptionA} holds.} There exists an invertible $\inR{V}{n}{n}$ such that
\begin{equation}
	\inv{V}S_\pi V=\diag(\varsigma_1,\dots,\varsigma_n), \label{thm:S:Vdef}
\end{equation}
where the columns of $V$ are eigenvectors of $S_\pi$ and the eigenvalues $\varsigma_i$ satisfy
\begin{align*}
\varsigma_i\in\left[\varsigma_\text{min},\varsigma_\text{max}\right]\eqdef\left[1,\,1+\frac{1}{\pi N(1+\pi)}\right].
\end{align*}
\end{theorem}
\begin{proof}
Pre- and post-multiplying~\eqref{eq:SandT} with $\trans{Z}$ and $Z$ from Theorem~\ref{thm:T} yields
\begin{align*}
\trans{Z}\tinv{R}S_\pi\trans{R}Z = \frac{1}{N-1}\left((1+\pi N)\trans{Z}T_{\pi}Z-\I\right).
\end{align*}
Since $\trans{Z}T_{\pi}Z=\diag(\tau_1,\dots,\tau_n)$, the matrix $\trans{Z}\tinv{R}S_\pi\trans{R}Z$ is diagonal. Set $V\eqdef\trans{R}Z$, which is invertible because $\det(R)\neq 0$ and $\xTx{Z}=\I$, then the columns of $V$ are eigenvectors of $S_{\pi}$ associated with eigenvalues $\varsigma_i=((1+\pi N)\tau_i-1)/(N-1)$. The bounds on $\varsigma_i$ are obtained from the bounds on $\tau_i$.
\end{proof}

The significance of Theorems~\ref{thm:T} and~\ref{thm:S} is that the diagonalizable matrices $T_{\pi}$ and $S_{\pi}$ have eigenvalues that are both bounded away from zero and contained in finite intervals, in contrast to the original formulation~\cite{HOGSVD} that requires a full rank condition and corresponds to $\pi =0$. More precisely, the range of eigenvalues of $S_\pi$ is contracted from $[1,\infty)$ for the original formulation to $\left[1,1+1/(\pi N(1+\pi))\right]$ in our case, which bounds the spectral condition number as $\kappa(S_\pi)\eqdef\twonorm{S_\pi}\twonorm{\inv{S_\pi}}\leq 1+1/(\pi N(1+\pi))$.

Before examining the eigenvalues of $S_{\pi}$ and $T_{\pi}$ further, we state our version of the HO-CSD and HO-GSVD. The HO-CSD and HO-GSVD have already been described in~\cite{LOANSLIDES} and~\cite{HOGSVD}, respectively, but our modified $D_{i,\pi}$ from~\eqref{eq:Di2} allows us to omit the requirements that $A_i$ and $Q_i$ be full rank.
\begin{definition}[HO-CSD]\label{def:HOCSD}
\added{Given $\,Q_1,\dots,\,Q_N$ satisfying~\eqref{eq:assumptionQ} and $N\geq 2$}, the HO-CSD of $\inR{Q_i}{m_i}{n}$ is given by $Q_i=U_i\Sigma_i\trans{Z}$, $i=1,\dots,N$, with ${Z}$ defined as in~\eqref{thm:T:Zdef}. The matrices $\inR{\Sigma_i}{n}{n}$ with $\Sigma_i=\diag(\sigma_{i,1},\dots,\sigma_{i,n})\succeq 0$ are obtained from
\begin{align*}
B_i \eqdef Q_iZ, \qquad B_i = \left[b_{i,1},\dots,b_{i,n}\right],\qquad \sigma_{i,k}=\twonorm{b_{i,k}},
\end{align*}
and $\inR{U_i}{m_i}{n}$ with $U_i=\left[u_{i,1},\dots,u_{i,n}\right]$ from
\[
u_{i,k} =
\begin{cases}
b_{i,k}/\sigma_{i,k} & \text{if } \sigma_{i,k}>0\\
u\in\R^{m_i}\text{ with }\twonorm{u}=1 & \text{if } \sigma_{i,k}=0.
\end{cases}
\]
\end{definition}
The left basis vectors $u_{i,k}$ have unit 2-norm and are, under certain circumstances, mutually orthogonal, \added{in which case they coincide with certain left generalized singular vectors of all pair-wise standard GSVD factorizations} (see Section~\ref{sec:commonsubspaces}). Because we allow for $\rank(Q_i)<n$, it is possible that $Q_iz_k = 0$ for some eigenvector $z_k$ of $T_\pi$,  consequently making the corresponding generalized singular value~$\sigma_{i,k} = 0$. In these cases, the column $u_{i,k}$ can be chosen freely or the corresponding row of~$\Sigma_i$ can be dropped. \added{Alternatively, they can be chosen to be orthogonal to all other columns, such as stated in the following lemma.}
\begin{newlemma}\label{thm:zeroSV}
Suppose that $\rank(Q_i)=n-K$ and let the generalized singular values $\sigma_{i,k}$ be ordered such that $\sigma_{i,1}=\dots=\sigma_{i,K}=0$ for $1\leq K < n$, and $\sigma_{i,j}>0$ for $j>K$. There exists $u_{i,1},\dots,u_{i,K}$ such that $\trans{u_{i,k}}u_{i,j}=0\,\forall k\leq K,\, j>K$.
\end{newlemma}
\ifArxivFormat
\begin{proof}
\else
\begin{newproof}
\fi
Let $\sigma_{i,1}=\dots=\sigma_{i,K}=0$ for $1\leq K < n$, and $\sigma_{i,j}>0$ for $j>K$. Then, $\spanv{u_{i,K+1},\dots,u_{i,n}}=\range(Q_i)$ and there exist $u_{i,1},\dots,u_{i,K}$ such that $\spanv{u_{i,1},\dots,u_{i,K}}=\ker(\trans{Q}_i)$, e.g. the columns of $\hat{U}_{i,2}$ associated with the standard SVD of $Q_i$:
\begin{align}\label{eq:zeroSV}
Q_i=\begin{bmatrix}\hat{U}_{i,1} & \hat{U}_{i,2}\end{bmatrix}\begin{bmatrix}\hat{\Sigma}_i & 0\\ 0 & 0\end{bmatrix}\begin{bmatrix}\hat{V}_{i,1} & \hat{V}_{i,2}\end{bmatrix}^\Tr,\qquad \hat{\Sigma}_i\succ 0.
\end{align}
\ifArxivFormat
\end{proof}
\else
\end{newproof}
\fi
\begin{definition}[HO-GSVD]\label{def:HOGSVD}
\added{Given $A_1,\dots,A_N$ satisfying~\eqref{eq:assumptionA} and $N\geq 2$}, the HO-GSVD of $\inR{A_i}{m_i}{n}$ is given by $A_i=U_i\Sigma_i\trans{V}$, with ${V}$ defined as in~\eqref{thm:S:Vdef}. The matrices $\inR{\Sigma_i}{n}{n}$ with $\Sigma_i=\diag(\sigma_{i,1},\dots,\sigma_{i,n})\succeq 0$ are obtained from
\begin{align*}
B_i \eqdef A_i\tinv{V}, \qquad B_i = \left[b_{i,1},\dots,b_{i,n}\right],\qquad \sigma_{i,k}=\twonorm{b_{i,k}},
\end{align*}
and $\inR{U_i}{m_i}{n}$ with $U_i=\left[u_{i,1},\dots,u_{i,n}\right]$ from
\[
u_{i,k} =
\begin{cases}
b_{i,k}/\sigma_{i,k} & \text{if } \sigma_{i,k}>0\\
u\in\R^{m_i}\text{ with }\twonorm{u}=1 & \text{if } \sigma_{i,k}=0.
\end{cases}
\]
\end{definition}
\added{%%%
According to Theorem~\ref{thm:S}, Definitions~\ref{def:HOCSD} and~\ref{def:HOGSVD} are equivalent in the sense that the HO-GSVD can be obtained from setting $V=\trans{R}Z$: 
\begin{align}\label{eq:GSVDCSD}
B_i = A_i\tinv{V} = Q_i R \inv{R} Z = Q_i Z,
\end{align}
where the rightmost term corresponds to $B_i$ as found in Definition~\ref{def:HOCSD}. The matrix of left basis vectors $U_i$ and the generalized singular values therefore depend only on the column space $Q$. However, when the HO-GSVD and the HO-CSD are computed separately, and $T_\pi$ and $S_\pi$ have eigenvalues with geometric multiplicity greater than $1$, it is possible that $V\neq \trans{R}Z$.
}%%%
\begin{newremark}
For rank-deficient $A_i$, the reader may wonder why the standard formulation of $S_\pi$ and $T_\pi$ with $\pi=0$ are not adapted by substituting the pseudoinverse for the inverse in~\eqref{eq:S} and~\eqref{eq:T}. The reason is that, in general, $\pinv{A_i}=\pinvbr{Q_i R}\neq \pinv{R}\pinv{Q_i}$ and using the pseudoinverse, the relationship~\eqref{eq:SandT} does not hold. However, relationship~\eqref{eq:SandT} is fundamental in determining the minimum and maximum eigenvalue of $S_\pi$ that will play an important role in subsequent sections, which is why pseudoinverses are not considered further.
\end{newremark}
\section{Common and isolated subspaces}\label{sec:commonsubspaces}%
\added{The HO-CSD and HO-GSVD identify directions, corresponding to columns of $Z$ and $V$, that, in the sense of~\eqref{eq:factorization2}, contribute equally to the factorizations of $Q_i$ and $A_i$, respectively. The directions are the right basis vectors $v_{i,k}$ associated with generalized singular values that are identical for each $Q_i$ and $A_i$, i.e. $\sigma_{i,k}=\sigma_{j,k}$}. These vectors form subspaces~\cite{HOGSVD, LOANSLIDES}, which are referred to as the common HO-CSD and HO-GSVD subspaces, and are defined in the following:
\begin{newdefinition}\label{def:commonsubspace}
The common HO-CSD subspace is defined as
\begin{align*}
\mathcal{T}_{N}\left\lbrace Q_1,\dots,Q_N\right\rbrace\eqdef \set{z\in\R^n}{T_{\pi}z=\tau_\text{min} z},
\end{align*}
and the common HO-GSVD subspace as
\begin{align*}
\mathcal{S}_{N}\left\lbrace A_1,\dots,A_N\right\rbrace\eqdef \set{v\in\R^n}{S_{\pi}v=\varsigma_\text{min} v},
\end{align*}
where $\tau_\text{min}$ and $\varsigma_\text{min}\!$ are the lower bounds on the range of eigenvalues defined in Theorems~\ref{thm:T} and~\ref{thm:S}, and $N\geq 2$.
\end{newdefinition}
Note that for a given set of matrices $A_1,\dots,A_N$, the subspaces $\mathcal{T}_{N}\left\lbrace Q_1,\dots,Q_N\right\rbrace$ and $\mathcal{S}_{N}\left\lbrace A_1,\dots,A_N\right\rbrace$ might be empty. \added{By Theorem~\ref{thm:S}, the HO-GSVD and HO-CSD subspaces are related by
\begin{align}\label{eq:HOCSDHOGSVDcommon}
\mathcal{S}_N&\left\lbrace A_1,\dots,A_N\right\rbrace = \set{\trans{R}z\in\R^n}{z\in\mathcal{T}_N\left\lbrace Q_1,\dots,Q_N\right\rbrace},
\end{align}
so that $v\in\mathcal{S}_{N}\left\lbrace A_1,\dots,A_N\right\rbrace$ iff $z=\tinv{R}v\in\mathcal{T}_{N}\left\lbrace Q_1,\dots,Q_N\right\rbrace$. The definition of the common subspace is complemented in the following theorem.}
\begin{newthm}\label{thm:HOCSDcommon}
The following statements are equivalent:
\begin{enumerate}[label={\thetheorem\alph*}]
\item $\mathcal{T}_{N}\left\lbrace Q_1,\dots,Q_N\right\rbrace\neq\emptyset$.\label{thm:HOCSDcommonI}
\item There exists $\hat{z}\in\R^n$ that is a standard right singular vector for each $Q_i$ and associated with a standard singular value $\hat{\sigma}=1/\sqrt{N}$ for each $Q_i$.\label{thm:HOCSDcommonII}
\item For each $Q_i$, there is a left basis vector $u_{i,k}$ satisfying $\trans{u_{i,k}}u_{i,p}=0$ $\forall p\neq k$ and the corresponding generalized singular values is $\sigma_{i,k}=1/\sqrt{N}$ for each $Q_i$.\label{thm:HOCSDcommonIII}
\end{enumerate}
\end{newthm}
\ifArxivFormat
\begin{proof}
\else
\begin{newproof}
\fi
The biconditional relationship \ref{thm:HOCSDcommonI}~$\Leftrightarrow$~\ref{thm:HOCSDcommonII} is a consequence of Theorem~\ref{thm:T}. Equality holds in~\eqref{eq:cauchyineq} iff $t$ is an eigenvector of $\xTx{Q_i}$~\cite[Thm. 7]{MEANS} or consequently in~\eqref{eq:ineqTa} iff $t$ is an eigenvector of each $\xTx{Q_i}$ for $i=1,\dots,N$. Equality holds in~\eqref{eq:ineqTb} iff $\trans{t}(\xTx{Q_i}+\pi I)t=\trans{t}(\xTx{Q_j}+\pi I)t$ for $i,j=1,\dots,N$. It follows that $T_{\pi}t=\tau_\text{min} t$ iff $t$ is a standard right singular vector for each $Q_i$ and from~\eqref{eq:assumptionQ} that $1=N\hat{\sigma}^2$, where $\hat{\sigma}=1/\sqrt{N}$ is the corresponding standard singular value. To show \ref{thm:HOCSDcommonII}~$\Rightarrow$~\ref{thm:HOCSDcommonIII}, let $\hat{u}_{i,k}$ be the corresponding standard left singular vector, then $Q_i z_k = \hat{\sigma}\hat{u}_{i,k}$ and from the HO-CSD, $Q_i z_k =\sigma_{i,k}u_{i,k}$, so the generalized singular values satisfy $\sigma_{i,k}=\hat{\sigma}$ since $\twonorm{\hat{u}_{i,k}}=\twonorm{u_{i,k}}=1$. To show that $\trans{u_{i,k}}u_{i,p}=0$ $\forall p\neq k$, consider the following equations for $\sigma_{i,p}\neq 0$:
\begin{align*}
\trans{u_{i,k}}u_{i,p} = \frac{\trans{b_{i,k}}b_{i,p}}{\sigma_{i,k}\sigma_{i,p}}
= \frac{\trans{z_{k}}\xTx{Q_i}z_{p}}{\sigma_{i,k}\sigma_{i,p}}
= \frac{\sigma_{i,k}}{\sigma_{i,p}}\trans{z_{k}}z_{p} = 0,
\end{align*}
where $b_{i,k}$ denotes column $k$ of the matrix $B_i$ from Definition~\ref{def:HOCSD}.

To show \ref{thm:HOCSDcommonIII}~$\Rightarrow$~\ref{thm:HOCSDcommonII}, suppose that \ref{thm:HOCSDcommonIII} holds and let $z_k$ be the corresponding right generalized singular vector. Then, $\xTx{Q_i}z_k = Z\Sigma_i U_i\trans{U_i}\Sigma_i Z z_k = \sigma_{i,k}^2 z_k$ since $\trans{u_{i,k}}u_{i,p}=0$ $\forall p\neq k$, hence $z_k$ is a shared standard right singular vector associated with a standard singular value $\sigma_{i,k}$.
\ifArxivFormat
\end{proof}
\else
\end{newproof}
\fi
\added{%
Note that statement~\ref{thm:HOCSDcommonIII} implies that the corresponding left basis vector $u_{i,k}$ is an eigenvector for $Q_i\trans{Q_i}$ for each $i$ and therefore also a \emph{standard} left singular vector for each $Q_i$.
}%

\added{The common HO-GSVD and HO-CSD subspaces are related by~\eqref{eq:HOCSDHOGSVDcommon}, and Theorem~\ref{thm:HOCSDcommon} can be adapted for the common HO-GSVD subspace as follows.}
\begin{newcorollary}\label{thm:HOGSVDcommon}
The following statements are equivalent:
\begin{enumerate}[label={\thetheorem\alph*}]
\item $\mathcal{S}_{N}\left\lbrace A_1,\dots,A_N\right\rbrace\neq\emptyset$.\label{thm:HOGSVDcommonI}
\item For each $A_i$, there is a left basis vector $u_{i,k}$ satisfying $\trans{u_{i,k}}u_{i,p}=0$ $\forall p\neq k$ and the corresponding generalized singular value is $\sigma_{i,k}=1/\sqrt{N}$ for each $A_i$.\label{thm:HOGSVDcommonIII}
\item There exists $v\in\R^n$ that is an eigenvector for each pairwise quotient $D_{i,\pi}\inv{D_{j,\pi}}$ associated with an eigenvalue $\lambda_{i,j}=1$.\label{thm:HOGSVDcommonIV}
\end{enumerate}
\end{newcorollary}
\ifArxivFormat
\begin{proof}
\else
\begin{newproof}
\fi
The biconditional relationship \ref{thm:HOGSVDcommonI}~$\Leftrightarrow$~\ref{thm:HOGSVDcommonIII} immediately follows from~\eqref{eq:HOCSDHOGSVDcommon} and Theorem~\ref{thm:HOCSDcommon}. To show~\ref{thm:HOGSVDcommonIII} $\Rightarrow$ \ref{thm:HOGSVDcommonIV}, substitute the HO-GSVD in~\eqref{eq:Di} to obtain
\begin{align}\label{eq:Didiag}
D_{i,\pi} = V\underbrace{\Sigma_i\xTx{U_i}\Sigma_i}_{\reqdef W_i}\trans{V}+\pi\xTx{A}
=V\left( W_i+\pi\sum_{p=1}^N W_p\right)\trans{V},
\end{align}
so that $D_{i,\pi}\inv{D_{j,\pi}}=V(W_i+\pi\sum_{p=1}^N W_p)\invbr{W_j+\pi\sum_{p=1}^N W_p}\inv{V}$. Because of~\ref{thm:HOGSVDcommonIII}, each $W_i$ has the block-diagonal form $W_i = \diag(\underline{W_i},\,\sigma_{i,k}^2+\pi\sum_{p=1}^N\sigma_{p,k}^2,\,\overline{W_i})$, where the scalar entry is on the $k$th row of $W_i$ and $\underline{W_i}$ and $\overline{W_i}$ are principal submatrices of $W_i$. Again from~\ref{thm:HOGSVDcommonIII}, $\sigma_{i,k}=\sigma_{j,k}$, so that $D_{i,\pi}\inv{D_{j,\pi}}v=v$. To complete the proof, we show that~\ref{thm:HOGSVDcommonIV}~$\Rightarrow$~\ref{thm:HOGSVDcommonI} by right-multiplying $S_\pi$ from~\eqref{eq:S} with $v$ from~\ref{thm:HOGSVDcommonIV}:
\begin{align*}
S_\pi v &= \frac{1}{N(N-1)}\sum_{i=1}^N\sum_{j=i+1}^N\left(\lambda_{i,j}v+\lambda_{j,i}v\right)
= \frac{2}{N(N-1)}\sum_{i=1}^N\sum_{j=i+1}^N v\\
&=\frac{2}{N(N-1)}\sum_{i=1}^N(N-i) v = \frac{2}{N(N-1)}\left(N^2-\frac{N^2+N}{2}\right) v = \varsigma_\text{min}v.
\end{align*}
\ifArxivFormat
\end{proof}
\else
\end{newproof}
\fi
\added{The ``common features'' of $A_i$ in~\eqref{eq:factorization2} can therefore be identified by the right basis vectors associated with eigenvalues of $S_\pi$ that equal $\varsigma_\text{min}$. In general, $\tinv{R}z$ is not an eigenvector for $\xTx{A_i}=R\xTx{Q_i}\trans{R}$, so that statement~\ref{thm:HOCSDcommonII} cannot be adapted to the HO-GSVD subspace, and the corresponding right basis vectors associated with the common subspace are not orthogonal in general. However, the right basis vectors spanning $\mathcal{S}_{N}\left\lbrace A_1,\dots,A_N\right\rbrace$ are eigenvectors of all pairwise quotients $D_{i,\pi}\inv{D_{j,\pi}}$, which is exploited in~\cite{LOANSLIDES} to compute the common HO-GSVD subspace using the standard pairwise GSVD. In addition, one can reformulate Statement~\ref{thm:HOGSVDcommonIV} to show that there exists a vector $\tilde{v}=\inv{D_{j,\pi}} v=\inv{D_{i,\pi}} v$, $v\in\mathcal{S}_{N}\left\lbrace A_1,\dots,A_N\right\rbrace$, that solves the \emph{higher-order generalized singular value problem} $\xTx{A_i}\tilde{v}=\mu\xTx{A_j}\tilde{v}$ with $\mu=1$.}
% \added{ This is exploited in the following lemma~\cite{LOANSLIDES}, which will be used in Section~\ref{sec:computation}.}
%\begin{newlemma}\label{thm:commonintersection}
%The common HO-CSD subspace can be characterized by
%\begin{align*}
%\mathcal{T}_{N}\left\lbrace Q_1,\dots,Q_N\right\rbrace =\bigcap_{i=2}^N\mathcal{S}_{2}\left\lbrace Q_{i-1},Q_i\right\rbrace.
%\end{align*}
%\end{newlemma}
%\begin{newproof}
%First, note that the QR factorization of the stacked $Q_1,\dots,Q_N$ yields $R=I$. By~\eqref{eq:HOCSDHOGSVDcommon}, it therefore holds that $\mathcal{T}_{N}\left\lbrace Q_1,\dots,Q_N\right\rbrace=\mathcal{S}_{N}\left\lbrace Q_1,\dots,Q_N\right\rbrace$. Next, by~\ref{thm:HOGSVDcommonIV} it holds that $\mathcal{S}_{N}\left\lbrace Q_1,\dots,Q_N\right\rbrace=\cap_{i=1}^N\cap_{j=i+1}^N\mathcal{S}_{2}\left\lbrace Q_i,Q_j\right\rbrace$. Representing $\mathcal{S}_{2}\left\lbrace Q_i,Q_j\right\rbrace$ as the intersection of two sets, e.g. $\mathcal{S}_{2}\left\lbrace Q_i,Q_j\right\rbrace=\mathcal{Q}_i\cap\mathcal{Q}_j$, and applying standard set intersection rules yields the desired result.
%\end{newproof}

In contrast to the common subspace, the isolated part of~\eqref{eq:factorization2} that is unique a single $A_i$ is identified by the right basis vectors associated with eigenvalues of $S_\pi$ ($T_\pi$) that equal $\varsigma_\text{max}$ ($\tau_\text{max}$).
\begin{newdefinition}\label{def:isolatedsubspace}
The isolated HO-CSD subspace is defined as
\begin{align*}
\mathcal{T}_{1}\left\lbrace Q_1,\dots,Q_N\right\rbrace\eqdef \set{z\in\R^n}{T_{\pi}z=\tau_\text{max} z},
\end{align*}
and the isolated HO-GSVD subspace as
\begin{align*}
\mathcal{S}_{1}\left\lbrace A_1,\dots,A_N\right\rbrace\eqdef \set{v\in\R^n}{S_{\pi}v=\varsigma_\text{max} v},
\end{align*}
where $\tau_\text{max}$ and $\varsigma_\text{max}\!$ are upper bounds on the range of eigenvalues defined in Theorems~\ref{thm:T} and~\ref{thm:S}, and $N\geq 2$.
\end{newdefinition}
\begin{newthm}\label{thm:HOCSDisolated}
The following statements are equivalent:
\begin{enumerate}[label={\thetheorem\alph*}]
\item $\mathcal{T}_{1}\left\lbrace Q_1,\dots,Q_N\right\rbrace\neq\emptyset$.\label{thm:HOCSDisolatedI}
\item There exists $\hat{z}\in\R^n$ that is a standard right singular vector for each $Q_i$ and associated with a standard singular value $\hat{\sigma}_{j,k}=1$ for one $Q_j$ and $\hat{\sigma}_{i,k}=0$ for all other $Q_i$, $i\neq j$.\label{thm:HOCSDisolatedII}
\item For each $Q_i$, there is a left basis vector $u_{i,k}$ satisfying $\trans{u_{i,k}}u_{i,p}=0$ $\forall p\neq k$, and the corresponding generalized singular value is $\sigma_{j,k}=1$ for one $Q_j$ and $\sigma_{i,k}=0$ for all other $Q_i$, $i\neq j$.\label{thm:HOCSDisolatedIII}
\end{enumerate}
\end{newthm}
\ifArxivFormat
\begin{proof}
\else
\begin{newproof}
\fi
The biconditional relationship \ref{thm:HOCSDisolatedI}~$\Leftrightarrow$~\ref{thm:HOCSDisolatedII} is a consequence of the proof of Theorem~\ref{thm:T}. According to Lemma~\ref{thm:P}, equality is attained in~\eqref{eq:ineqTupper} iff for each summand, $t\in\spanv{v_1^i,\dots,v_p^i}$, where $v_k^i$ are eigenvectors of $\xTx{Q_i}$ associated with eigenvalues that are equal to either $0$ or $1$. It remains to consider~\eqref{eq:assumptionQ}. The proof for~\ref{thm:HOCSDisolatedII}~$\Leftrightarrow$~\ref{thm:HOCSDisolatedIII} follows the proof of Theorem~\ref{thm:HOCSDcommon}. For the case $\hat{\sigma}_{i,k}=\sigma_{i,k}=0$, it is assumed Lemma~\ref{thm:zeroSV} has been applied.
\ifArxivFormat
\end{proof}
\else
\end{newproof}
\fi
\added{By Theorem~\ref{thm:S}, the isolated HO-GSVD and HO-CSD subspaces are related by}
\begin{align}\label{eq:HOCSDHOGSVDisolated}
\mathcal{S}_1&\left\lbrace A_1,\dots,A_N\right\rbrace = \set{\trans{R}z\in\R^n}{z\in\mathcal{T}_1\left\lbrace Q_1,\dots,Q_N\right\rbrace},
\end{align}
and Theorem~\ref{thm:HOCSDisolated} is reformulated for the HO-GSVD as follows.
\begin{newcorollary}\label{thm:HOGSVDisolated}
The following statements are equivalent:
\begin{enumerate}[label={\thetheorem\alph*}] % NEED TO UPDATE THIS
\item $\mathcal{S}_{1}\left\lbrace A_1,\dots,A_N\right\rbrace\neq\emptyset$.\label{thm:HOGSVDisolatedI}
\item For each $A_i$, there is a left basis vector $u_{i,k}$ satisfying $\trans{u_{i,k}}u_{i,p}=0$ $\forall p\neq k$, and the corresponding generalized singular value is $\sigma_{j,k}=1$ for one $A_j$ and $\sigma_{i,k}=0$ for all other $A_i$, $i\neq j$.\label{thm:HOGSVDisolatedIII}
\item There exist $v\in\R^n$ and $i\in\lbrace 1,\dots,N\rbrace$ such that $v$ is an eigenvector for each pairwise quotient $D_{p,\pi}\inv{D_{j,\pi}}$ associated with eigenvalues $\lambda_{i,j}=\frac{1+\pi}{\pi}$, $\lambda_{j,i}=\frac{\pi}{1+\pi}$ and $\lambda_{p,j}=\lambda_{j,p}=1$ for $j=\lbrace 1,\dots,N\rbrace$, $p=\lbrace 1,\dots,N\rbrace$ and $j\neq p\neq i$.\label{thm:HOGSVDisolatedIV}
\end{enumerate}
\end{newcorollary}
\ifArxivFormat
\begin{proof}
\else
\begin{newproof}
\fi
The proof follows the proof of Corollary~\ref{thm:HOGSVDcommon}. To show~\ref{thm:HOGSVDisolatedI}~$\Leftrightarrow$~\ref{thm:HOGSVDisolatedIII}, use~\eqref{eq:HOCSDHOGSVDisolated} and apply Theorem~\ref{thm:HOCSDisolated}. To show~\ref{thm:HOGSVDisolatedIII}~$\Rightarrow$~\ref{thm:HOGSVDisolatedIV}, use~\eqref{eq:Didiag}. Finally, to show~\ref{thm:HOGSVDisolatedIV}~$\Rightarrow$~\ref{thm:HOGSVDisolatedI}, compute $S_\pi v$ and assume without loss of generality that $i=1$:
\begin{align*}
S_\pi v &= \frac{1}{N(N-1)}\sum_{j=2}^N\left(\lambda_{1,j}+\lambda_{j,1}\right)v+\frac{1}{N(N-1)}\sum_{p=2}^N\sum_{j=p+1}^N\left(\lambda_{p,j}+\lambda_{j,p}\right)v\\
&= \frac{1}{N(N-1)}\sum_{j=2}^N\left(\frac{1+\pi}{\pi}+\frac{\pi}{1+\pi}\right)v+\frac{1}{N(N-1)}\sum_{p=2}^N\sum_{j=p+1}^N2v\\
&=\frac{1}{N}\left(\frac{1+\pi}{\pi}+\frac{\pi}{1+\pi}+N-2\right)v=\varsigma_\text{max}v.
\end{align*}
\ifArxivFormat
\end{proof}
\else
\end{newproof}
\fi
%\begin{newlemma}\label{thm:isolatedintersection}
%The isolated HO-CSD subspace can be characterized by
%\begin{align*}
%\mathcal{T}_{1}\left\lbrace Q_1,\dots,Q_N\right\rbrace =\bigcap_{i=2}^N\mathcal{S}_{1}\left\lbrace Q_{i-1},Q_i\right\rbrace.
%\end{align*}
%\end{newlemma}
%\begin{newproof}
%Adapt the proof of Lemma~\ref{thm:commonintersection} by substituting~\eqref{eq:HOCSDHOGSVDisolated} for~\eqref{eq:HOCSDHOGSVDcommon} and statement~\ref{thm:HOGSVDisolatedIV} for~\ref{thm:HOGSVDisolatedIV} to obtain the desired result.
%\end{newproof}
\added{Statements~\ref{thm:HOCSDcommonIII} and~\ref{thm:HOCSDisolatedIII} of Theorems~\ref{thm:HOCSDcommon} and~\ref{thm:HOCSDisolated} show that, in certain cases, the orthogonality of the left factor matrix, which always holds for the standard SVD and GSVD, is preserved for higher-order datasets (see also Section~\ref{sec:comparison}). If the generalized singular values $\sigma_{i,k}$, the left basis vectors $u_{i,k}$ and the right basis vectors $v_{k}$ are grouped according to whether they are associated with the common subspace ($k\in\mathcal{I}_N$), the isolated subspace ($k\in\mathcal{I}_1$) or neither of the subspaces ($k\in\mathcal{I}_\perp$), Definition~\ref{def:HOGSVD} can be refined as
\begin{align}\label{eq:HOGSVD2}
A_i = \begin{bmatrix}
U_{i,\mathcal{I}_1} & U_{i,\mathcal{I}_\perp} & U_{i,\mathcal{I}_N}
\end{bmatrix}
\begin{bmatrix}
\Sigma_{i,\mathcal{I}_1} \\ & \Sigma_{i,\mathcal{I}_\perp}\\ & & I/\sqrt{N}
\end{bmatrix}
\trans{\begin{bmatrix}
V_{\mathcal{I}_1} & V_{\mathcal{I}_\perp} & V_{\mathcal{I}_N}
\end{bmatrix}},
\end{align}
where $\Sigma_{i,\mathcal{I}_1}$ contains the generalized singular values associated with $\mathcal{S}_1\left\lbrace A_1,\dots,A_N\right\rbrace$ and $\Sigma_{i,\mathcal{I}_\perp}\succ 0$. In the notation of~\eqref{eq:HOGSVD2}, the three blocks of left basis vectors are mutually orthogonal, e.g. $\transbr{U_{i,\mathcal{I}_N}}U_{i,\mathcal{I}_1} =0$, which follows from statements~\ref{thm:HOGSVDcommonIII} and~\ref{thm:HOGSVDisolatedIII} of Corollaries~\ref{thm:HOGSVDcommon} and~\ref{thm:HOGSVDisolated}. Note that for the HO-GSVD, the right basis vectors are \emph{not} orthogonal in general.
}%

As can also be concluded from Theorems~\ref{thm:HOCSDcommon} and~\ref{thm:HOCSDisolated}, the parameter $\pi$ does not alter the common and isolated subspaces, which shows that the standard HO-GSVD formulation and the present one are equivalent.
\begin{corollary}\label{thm:samesub}
The common and isolated HO-GSVD and HO-CSD subspaces are independent of the value of $\pi$.
\end{corollary}
\begin{proof}
For the HO-CSD, the claim follows from statements \ref{thm:HOCSDcommonII} and~\ref{thm:HOCSDisolatedII} of Theorems~\ref{thm:HOCSDcommon} and~\ref{thm:HOCSDisolated}, which are independent of the value of $\pi$. \added{As a consequence of~\eqref{eq:HOCSDHOGSVDcommon} and~\eqref{eq:HOCSDHOGSVDisolated}, the claim is also true for the HO-GSVD.}
\end{proof}
\added{Note that Corollary~\ref{thm:samesub} ignores potential numerical inaccuracies, which are treated in Section~\ref{sec:computation}. Numerical inaccuracies can also cause rank deficiencies of the stacked matrix $A$, and the following Remark~\ref{rem:rankdeficient} explains how the HO-GSVD can be applied even when $A$ does not satisfy~\eqref{eq:assumptionA}.}
\begin{newremark}\label{rem:rankdeficient}
Suppose that assumption~\eqref{eq:assumptionA} does \emph{not} hold and that $\rank(A)=r<n$. Then, $S_\pi$ is undefined and~\eqref{eq:SandT} invalid. Let $\spanv{v_1,\dots,v_{n-r}}=\ker(A)$ be an orthogonal basis and set $A_{N+1}\eqdef\left[ v_1,\, \dots,\, v_{n-r}\right]^\Tr$. The HO-GSVD can be applied to the augmented dataset $A_1,\dots,A_{N+1}$, and at least $n-r$ directions of the resulting isolated HO-GSVD subspace are associated with $\ker(A)$.
\end{newremark}%
\section{The parameter $\pi$}\label{sec:pi}
The eigenvectors of $T_\pi$ that are in the common or isolated HO-CSD subspaces are not affected by the choice of $\pi$, but other (normalized) eigenvectors can be modified as $\pi$ varies. Here, we are interested in the limits of these eigenvectors as $\pi\rightarrow 0$ and $\pi\rightarrow \infty$.
%
%From~\eqref{eq:T}, it holds that for $\pi\rightarrow \infty$, $\anynorm{S_\pi-I}\rightarrow 0$ and $\anynorm{T_\pi}\rightarrow 0$, in which case any set of $n$ linearly independent vectors are eigenvectors for $S_\pi$ and $T_\pi$ at their limit points.
Since from~\eqref{eq:T} it holds that \added{$\lim_{\pi\rightarrow \infty} S_\pi = I$ and  $\lim_{\pi\rightarrow \infty} T_\pi = 0$}, some caution is required in determining the limits of the associated eigenvectors.

Semisimple eigenvalues are expected to be associated with the common or isolated subspaces and therefore not considered further. To examine the remaining eigenvectors associated with simple eigenvalues, we will make use of the following result in both cases:

\begin{theorem}[{\cite[Thm. 7 \& 8, Ch. 9.3, p. 130]{LAX}}]\label{thm:laxsimpel}
Let $M(x)$ be a differentiable square matrix-valued function of the real variable $x$. Suppose that $M(0)$ has a simple eigenvalue $m_0$. Then for $x$ small enough, $M(x)$ has an eigenvalue $m(x)$ that depends differentiably on $x$ with $m(0) = m_0$ and we can choose an eigenvector $h(x)$ of $M(x)$ pertaining to the eigenvalue $m(x)$ that depends differentiably on $x$.
\end{theorem}

\textbf{The case $\pi \to \infty$:}

\begin{lemma}[Eigenvectors of $T_{\pi}$ for $\pi\rightarrow\infty$]\label{thm:inflimit}
Consider the matrix $\tilde{T}_\infty$,
\begin{align}\label{eqn:Ttildeinf}
\tilde{T}_\infty \eqdef \frac{1}{N}\sum_{i=1}^N (\xTx{Q_i})^2,
\end{align}
and suppose that $\tilde{T}_\infty$ has a simple eigenvalue $\tilde{\tau}_\infty$ associated with an eigenvector $\tilde{z}_\infty$. Then, there exists an eigenvector $z(\pi)$ of $T_\pi$ that depends differentiably on $\pi$ and converges to $\tilde{z}_\infty$ as $\pi\rightarrow\infty$.
\end{lemma}
\begin{proof}
Use the Neumann series $\invbr{I-M} = \sum_{k=0}^{\infty} M^k$ with $\anynorm{M}<1$~\cite[Ch. 1.4]{KASO} \added{to expand each of the summands in~\eqref{eq:T} as $\invbr{\xTx{Q_i}+\pi I}=\frac{1}{\pi}\sum_{k=0}^{\infty}\left(\frac{-1}{\pi}\xTx{Q_i}\right)^k$, and rewrite $T_\pi$ as}
\begin{align*}
T_\pi=\frac{1}{N\pi}\sum_{i=1}^N\sum_{k=0}^{\infty}\left(\frac{-1}{\pi}\xTx{Q_i}\right)^k= \frac{1}{\pi} I - \frac{1}{N\pi^2}I+\sum_{i=1}^N\frac{1}{N\pi^3}(\xTx{Q_i})^2+\co{\frac{1}{\pi^4}},
\end{align*}
where $\anynorm{\frac{1}{\pi}\xTx{Q_i}}<1$ for $\pi>1$.  Set $\tilde{T}(\pi) =  \pi^3\, (T_\pi-\frac{N\pi-1}{N\pi^2}I)$, which depends differentiably on $\pi$ for $\pi>0$ and has the same eigenvectors as $T_\pi$. Neglecting higher-order terms $\co{1/\pi^4}$, the limit $\lim_{\pi\to\infty}\tilde{T}(\pi) = \tilde{T}_\infty$ is obtained, where equality holds element-wise.
Finally, defining
\begin{align*}
M(x)=\begin{cases}\tilde{T}_\infty & x=0,\\ \tilde{T}(1/x), & x> 0,\end{cases}
\end{align*}
the proof follows from Theorem~\ref{thm:laxsimpel}.
\end{proof}

\added{According to Lemma~\ref{thm:inflimit}, the eigenvectors of $T_\pi$ associated with simple eigenvalues can be chosen such that they converge to those of $\tilde{T}_\infty$ for large $\pi$. Suppose that some $Q_i$ has a ``dominant'' standard right singular vector $\bar{v}$, in the sense that $\trans{\bar{v}}\xTx{Q_i}\bar{v}\gg\trans{\bar{v}}\xTx{Q_j}{\bar{v}}$ for $j\neq i$. In this case, $\tilde{T}_\infty$ can be rewritten as $\tilde{T}_\infty=\xTx{Q_i}/N+\Delta$ with $\twonorm{\Delta}\ll\twonorm{\xTx{Q_i}/N}$. According to standard perturbation theory~\cite[Ch. 7.2.5]{GOLUB4}, $T_\pi$ will have an eigenvector $v=\bar{v}+\delta v$ with $\twonorm{\delta v}\ll 1.$} By using the orthogonality property $\sum_{i=1}^N \xTx{Q_i}=I$, the matrix $\tilde{T}_\infty$ defined in \eqref{eqn:Ttildeinf} can be rewritten as
\begin{equation}\label{eq:symprod}
\begin{aligned}
\tilde{T}_\infty
&=\frac{1}{N}\sum_{i=1}^{N}\xTx{Q_i}(I-\sum_{j\neq i}\xTx{Q_j})\\
&=\frac{1}{N}I -\frac{1}{N}\sum_{i=1}^{N-1}\sum_{j=i+1}^{N}\left(\xTx{Q_i}\xTx{Q_j}+\xTx{Q_j}\xTx{Q_i}\right),
\end{aligned}\end{equation}
where the summands on the second line are referred to as \emph{symmetrized products} or \emph{Jordan products} of $\xTx{Q_i}$ and $\xTx{Q_j}$~\cite[Ch. 10]{LAX}. \added{The form~\eqref{eq:symprod} shows that the eigenvectors of $T_\pi$ will also converge to those of a ``dominant'' symmetrized product for large $ \pi$.}

\textbf{The case $\pi \to 0$:}

For any rank-deficient $Q_i$ and $\pi=0$, the corresponding term $\xTx{Q_i}+\pi I$ appearing in the definition of $T_\pi$ in \eqref{eq:T} is singular. However, by using the standard SVD $\xTx{Q_i}=V_i\diag(\sigma_{i,1}^2,\dots,\sigma_{i,r}^2,0,\dots,0)\trans{V_i}$ with $r=\rank(Q_i)$, one can show that
\begin{align*}
\lim_{\pi\rightarrow 0}\pi\invbr{\xTx{Q_i}+\pi I}=V_i\diag(\underbrace{0,\dots,0}_{r\text{ times}},\underbrace{1,\dots,1}_{n-r\text{ times}})\trans{V_i},
\end{align*}
where this limit is zero if $Q_i$ is instead full rank.
The following lemma provides useful information in the  case where some of the $Q_i$ are rank-deficient.
\begin{lemma}[Eigenvectors of $T_{\pi}$ for $\pi\rightarrow 0$]\label{thm:zerolimit}
Suppose that some of the $Q_i$ are rank-deficient. Consider
\begin{align}\label{eqn:Ttilde0}
\tilde{T}_0 \eqdef \frac{1}{N}\sum_{i=1}^N Q_i^\dagger Q_i,
\end{align}
where $Q_i^\dagger=\lim_{\pi\rightarrow 0} \trans{Q_i}\invbr{\pi I + Q_i\trans{Q_i}}$ is the Moore-Penrose pseudoinverse of $Q_i$~\cite[P5.5.2]{GOLUB4},
and suppose that $\tilde{T}_0$ has a simple eigenvalue $\tilde{\tau}_0$ associated with an eigenvector $\tilde{z}_0$. Then, there exists an eigenvector $z(\pi)$ of $T_\pi$ that depends differentiably on $\pi$ and converges to $\tilde{z}_0$ as $\pi\rightarrow 0$.
\end{lemma}
\begin{proof}
Use the Woodbury matrix identity~\cite[Ch. 2.1.4]{GOLUB4} to rewrite $\pi T_\pi$ for $\pi>0$ as
\begin{align}\label{eqn:piTpi}
\pi T_\pi &= \frac{1}{N}\sum_{i=1}^N \left(I-\trans{Q_i}\invbr{\pi I+Q_i\trans{Q}_i} Q_i\right),
\end{align}
with $\lim_{\pi\rightarrow 0} \pi\, (T_\pi-\frac{1}{\pi}I) = -\tilde{T}_0$ (element-wise), where $\tilde{T}_0$ and $\pi\, (T_\pi-\frac{1}{\pi}I)$ share the same eigenspace~\cite[Ch.2]{GOLUB4}. Differentiability of the matrix $\pi T_\pi$ with respect to $\pi$ at $0$ is easily shown by substitution of the standard SVD of each $Q_i$ into \eqref{eqn:piTpi}. The proof then follows from Theorem~\ref{thm:laxsimpel}.
% \begin{align*}
% \frac{d}{d\pi}(\pi T_\pi)=\frac{1}{N}\sum_{i=1}^N\invbr{\xTx{Q_i}+\pi I}\xTx{Q_i}\invbr{\xTx{Q_i}+\pi I}
% \end{align*}
% where we used $\frac{d}{dx}\inv{M}(x)=-\inv{M}(x)(\frac{d}{dx} M(x))\inv{M}(x)$~\cite[Ch. 9.1]{LAX}. It holds that $\frac{d}{d\pi}(\pi T_\pi)= \frac{1}{N}\sum_{i=1}^N Q_i^\dagger(\trans{Q_i})^\dagger$ for $\pi\rightarrow 0$ and is therefore applicable.
\end{proof}

Note that if \emph{all} $Q_i$ have full \added{column} rank, then $Q_i^\dagger Q_i = I$ and $\tilde{T}_0$ has no simple eigenvalues.
The matrix $Q_i^\dagger Q_i$ is the orthogonal projector onto $\range(\trans{Q_i})$ and $Q_i^\dagger Q_i=I$ if $Q_i$ has full rank. It follows that if some $Q_j$ are rank deficient, then the eigendecomposition of $T_\pi$ can be chosen such that it equals the eigendecomposition of the sum of projectors onto $\range(\trans{Q_j})$ (the orthogonal complement of $\kernel(Q_j)$), but if all $Q_i$ have full rank, then the eigenvectors of $\lim_{\pi\rightarrow 0} T_\pi$ are those of $T_0$, i.e.~\eqref{eq:T} with $\pi=0$.

The limits for $S_\pi$ can be obtained from pre- and post-multiplying $\tilde{T}_0$ or $\tilde{T}_\infty$ with $\trans{R}$ and $\tinv{R}$, respectively.
%\begin{lemma}[Ordering of $T_\pi$]
%If $\pi_2 > \pi_1 > 0$, then $T_{\pi_1}\succ T_{\pi_2}$.
%\end{lemma}
%\begin{proof}
%Consider the summands of $T_\pi$ and define $F\eqdef\xTx{Q_i}+\pi_1 I$, $G\eqdef\xTx{Q_i}+\pi_2 I$ and $H\eqdef G-F=(\pi_2-\pi_1)I\succ 0$. It holds that $\inv{F}=\invbr{G-H}=\invbr{(I-H\inv{G})G}=\inv{G}\invbr{I-H\inv{G}}$, which, after noting that $\anynorm{H\inv{G}}=(\pi_2-\pi_1)/\anynorm{G}< 1$ for sufficiently small $\pi_2-\pi_1$, can be expanded using a Neumann series as
%\[
%\inv{F}=\inv{G}\sum_{j=0}^\infty \left(H\inv{G}\right)^n =\inv{G}+\sum_{j=1}^\infty \inv{G}\left(H\inv{G}\right)^j.
%\]
%From $\inv{G}\left(H\inv{G}\right)^j\succ 0$, it follows that $\inv{F}-\inv{G}=\invbr{\xTx{Q_i}+\pi_1 I}-\invbr{\xTx{Q_i}+\pi_2 I}\succ 0$ and hence $T_{\pi_1}\succ T_{\pi_2}$.
%\end{proof}

Apart from rotating the eigenvectors, the choice of $\pi$ also affects the function $f_\pi:\R^{n}\rightarrow \Rpos$,
\begin{align}\label{eq:f}
f_\pi(v) &= \frac{1}{N(N-1)}\sum_{i=1}^{N-1}\sum_{j=1}^N\left(
\frac{\trans{v}(\xTx{A_i}+\pi\xTx{A})v}{\trans{v}(\xTx{A_j}+\pi\xTx{A})v}+\frac{\trans{v}(\xTx{A_j}+\pi\xTx{A})v}{\trans{v}(\xTx{A_i}+\pi\xTx{A})v}\right),
\end{align}
where $\twonorm{v}=1$ and $f_\pi(v)\geq 1$. The function $f_\pi(v)$ measures the arithmetic mean of amplifications in a particular direction $v$ and has been shown to be related to the (HO-)GSVD~\cite{GSVDVAR2, GSVDVAR, LOANSLIDES}. For $N=2$, $\pi=0$ and full-rank $A_1$ and $A_2$, the gradient is zero for vectors that lie in the common HO-GSVD subspace~\cite{LOANSLIDES}, which can be extended to the isolated HO-GSVD subspace (Appendix~\ref{app:rayleigh}). The parameter $\pi$ has the effect of \emph{flattening out} $f_\pi$ and, in particular, \added{removing the singularities of $f_\pi(v)$ associated with the nullspace of $A_i$ for $\pi>0$ and $\rank(A)=n$, since in that case $\trans{v}(\xTx{A_i}+\pi\xTx{A})v > 0$ for $v\neq 0$.}

\section{Comparison with standard HO-GSVD, GSVD and SVD}\label{sec:comparison}
When one out of two matrices is the identity matrix, the GSVD reduces to the standard SVD~\cite{LOAN_GSVD}. The same has been shown for the full-rank HO-GSVD~\cite{HOGSVD}. When $N-1$ matrices $A_i$ are identity matrices, then the full-rank HO-GSVD reverts to the standard SVD of $A_j$, $j\neq i$. Here, this fact is demonstrated for our HO-GSVD as given in Definition~\ref{def:HOGSVD}.
\begin{theorem}\label{thm:HOGSVD_SVD}
Let $A_1$ be an arbitrary matrix and \added{$A_2=\dots=A_N=I$ with $N\geq 2$}. The HO-GSVD of $A_1,A_2,\dots,A_N$ with \added{$\pi>0$} yields the standard SVD of $A_1$.
\end{theorem}
\ifArxivFormat
\begin{proof}
\else
\begin{newproof}
\fi
Substitute the standard SVD $\hat{U}_1\hat{\Sigma}_1\trans{\hat{V}}_1=A_1$ and $A_j=I$, $j=2,\dots,N$, in~\eqref{eq:Di}, so that 
\begin{align*}
\trans{\hat{V}_1}D_1\hat{V}_1=(1+\pi)\xTx{\hat{\Sigma}_1}+\pi(N-1)I,\quad
\trans{\hat{V}_1}D_j\hat{V}_1=\pi\xTx{\hat{\Sigma}_1}+(1+\pi(N-1))I.
\end{align*}
The summands $D_{i,\pi}\inv{D_{j,\pi}}+D_{j,\pi}\inv{D_{i,\pi}}$ in the definition of $S_\pi$~\eqref{eq:S} are therefore diagonalized by $\hat{V}_1$, and $V=\hat{V}_1$ is an orthogonal eigenbasis for $S_\pi$. According to Definition~\ref{def:HOGSVD}, the HO-GSVD $A_1=U_1\Sigma_1\trans{V}$ is obtained from $B_1=A_1\tinv{V}=A_1\bar{V}_1=\bar{U}_1\bar{\Sigma}_1$, so that $U_1=\hat{U}_1$ and $\Sigma_1=\hat{\Sigma}_1$.
\ifArxivFormat
\end{proof}
\else
\end{newproof}
\fi
The HO-GSVD from Definition~\ref{def:HOGSVD} can also be related to the GSVD. For the special case that $N=2$, $\inR{A_1}{m_1}{n}$ with $m_1\geq n$ and $\rank(A_1) =n$ and an arbitrary $\inR{A_2}{m_2}{n}$, it can be shown that the HO-GSVD yields $\Sigma_i$ with $\xTx{\Sigma_1}+\xTx{\Sigma_2}=\I$ and orthogonal $U_1$ and $U_2$.
\begin{theorem}\label{thm:HOGSVD_GSVD2}
For $N=2$ and $\pi>0$, the HO-CSD from Definition~\ref{def:HOCSD} yields the standard CSD and the HO-GSVD from Definition~\ref{def:HOGSVD} yields the standard GSVD.
\end{theorem}
\ifArxivFormat
\begin{proof}
\else
\begin{newproof}
\fi
Since $\invbr{\xTx{Q_i}+\pi I\,}$ and  $\xTx{Q_i}$ with $i=1,2$ and $\xTx{Q_1}+\xTx{Q_2}=\I$ share the same eigenspace for any $\pi \in\Rpos$~\cite[Ch. 2]{GOLUB4}, the eigenvectors $z_k$ for $T_{\pi}$ can be chosen such that they are right singular vectors for $Q_1$ and $Q_2$. Let $b_{i,k}$ denote the columns of $B_i=Q_iZ$, then for $j\neq k$,
$\trans{b_{i,k}}b_{i,j}=\trans{z_k}\xTx{Q_i}z_j=\hat{\sigma}_{i,j}^2\trans{\bar{u}_{i,k}}\hat{u}_{i,j}=0$,
where $\hat{\sigma}_\times$ and $\hat{u}_\times$ denote standard singular values and left singular vectors, respectively. Hence, from $U_i\Sigma_i=B_i$, the columns of $U_i$ are either zero or orthonormal. Substituting $Q_i=U_i\Sigma_i\trans{V}$ in $\xTx{Q_1}+\xTx{Q_2}=\I$ yields
$
Z\xTx{\Sigma_1}\trans{Z}+Z\xTx{\Sigma_2}\trans{Z}=\I,
$
and from $\xTx{Z}=\I$, follows $\xTx{\Sigma_1}+\xTx{\Sigma_2}=\I$. The claim on the HO-GSVD follows from Theorem~\ref{thm:HOGSVD_GSVD2} with $V=\trans{R}Z$.
\ifArxivFormat
\end{proof}
\else
\end{newproof}
\fi
\begin{remark} Lemma~\ref{thm:HOGSVD_GSVD2} shows that for $N=2$ the three matrices, $T_{\pi}$, $\xTx{Q_1}$ and  $\xTx{Q_2}$, share the same eigenspace, but not every eigendecomposition of $T_{\pi}$ yields eigenvectors that are parallel to those of $\xTx{Q_1}$ and  $\xTx{Q_2}$. For example, suppose that $\dim(\ker(Q_i))=1$ and that $q_i\in\ker(Q_i)$, $i=1,2$, are linearly independent. \added{From pre- and post-multiplying $\xTx{Q_1}+\xTx{Q_2}=\I$ with $\trans{q}_1$ and $q_2$}, it holds that $\trans{q_1}q_2=0$. It follows that $\dim(\mathcal{T}_1)=2$, so that $T_{\pi}$ has a semisimple eigenvalue. When the associated eigenvectors are computed using numerical software, these will not necessarily be parallel to $q_1$ and $q_2$, and the HO-CSD will not necessarily yield orthonormal matrices $U_i$.
\end{remark}
The HO-GSVD from Definition~\ref{def:HOGSVD} can also be compared with the full-rank HO-GSVD~\cite{HOGSVD}. For $N=2$ and full-rank matrices $A_i$, both HO-GSVDs have been shown to be equivalent to the GSVD. Both HO-GSVDs have also been shown to yield the SVD of $A_j$ when $A_i=\I$ for $i\neq j$. For $N>2$, however, the HO-GSVD from Definition~\ref{def:HOGSVD} and~\cite{HOGSVD} will in general \emph{not} yield identical factorizations $A_i=U_i\Sigma_i\trans{V}$, even when $\rank (A_i)=n$. This can be seen by comparing the eigenspaces of $T_{\pi}$ from~\eqref{eq:T} for varying $\pi$, where $\pi =0$ corresponds to the standard HO-CSD~\cite{LOANSLIDES}. For $N=2$, the eigenvectors of $T_{\pi}$ are independent of the value of $\pi $ because its eigenvectors are fixed by the orthogonality property $\xTx{Q_1}+\xTx{Q_2}=\I$, while for $N>2$ this property is lost. From Theorem~\ref{thm:S}, it follows that the same holds for the HO-GSVD. However, it can be shown that in case the matrices $A_i$ and $Q_i$ have full rank, then the common HO-CSD and HO-GSVD subspaces will be the same for any value of $\pi$ (Corollary~\ref{thm:samesub}) \added{and $N>2$. Moreover, it follows from Theorem~\ref{thm:laxsimpel} that in the full-rank case, the eigenvectors of $T_\pi$ converge to those of the standard HO-GSVD as $\pi\rightarrow 0$.}

\section{\added{Computing the HO-GSVD}\label{sec:computation}}
The early literature on the standard GSVD ($N=2$) identified numerical issues for the case that $A$ from~\eqref{eq:A} and therefore $R$ are ill-conditioned~\cite{STEWART_EV,LOAN_GSVD,PAIGE_GSVD}. This problem was resolved by basing the GSVD computation on the CSD, hereby avoiding computing the inverse of $R$. To compute the full HO-GSVD~\eqref{eq:factorization}, we propose to use Algorithm~\ref{alg:HOGSVD}, which is based on the HO-CSD. An experimental Matlab implementation is provided in~\cite{HOGSVDGITHUB}.
\begin{algorithm}
 \caption{HO-GSVD Computation}\label{alg:HOGSVD}
 \begin{algorithmic}[1]
    \REQUIRE $A_1,\dots,A_N$, $\pi > 0$
    \ENSURE Factorizations $A_i=U_i\Sigma_i \trans{V},\, i=1,\dots,N$
    \STATE Obtain $Q_i R=A_i$ for $i=1,\dots,N$ from~\eqref{eq:A}\COMMENT{2Mn^2}
	\STATE Form $T_\pi$ using~\eqref{eq:T}\COMMENT{M n^2 + N n^3}
	\STATE Obtain the eigenvectors $z_1,\dots,z_n$ of $T_\pi$\COMMENT{n^3}\label{alg:HOGSVD:EV}
	\STATE Determine $\mathcal{I}_1$ and align $z_k$, $k\in\mathcal{I}_1$\COMMENT{2Mn^2+n^3}
	\FOR{$i = 1,\dots,N$ and $k = 1,\dots,n$}\label{alg:HOGSVD:loop}
  	  	\IF{$\sigma_{i,k}=\twonorm{Q_i z_k}>0$}
  	  		\STATE Set $u_{i,k} = Q_i z_k / \sigma_{i,k}$
  	  	\ELSE
  	  		\STATE Assign a column of $\hat{U}_{i,2}$ from~\eqref{eq:zeroSV} to $u_{i,k}$
  	  	\ENDIF
	\ENDFOR\COMMENT{Mn^2}
	\STATE Set $V=\trans{R}\begin{bmatrix}z_{1},\dots,z_{n}\end{bmatrix}$\COMMENT{n^3}\label{alg:HOGSVD:V}
 \end{algorithmic}
\end{algorithm}

Given a dataset $A_1,\dots,A_N$ and a parameter $\pi > 0$, Algorithm~\ref{alg:HOGSVD} starts by computing the thin QR factorization~\eqref{eq:A}, which enables use of the HO-CSD to avoid computing the inverse of a potentially ill-conditioned $R$. Next, the terms $\invbr{\xTx{Q_i}+\pi I}$ are computed to obtain $T_\pi$. Forming the products $\xTx{Q_i}$ may lead to a loss of accuracy if $Q_i$ has small singular values, but the condition number $\kappa(\xTx{Q_i}+\pi I)$ can be controlled by choosing $\pi$ as follows. Let $\kappa_\text{max} > 1$ and $\hat{\sigma}_{i,\text{min}}$ and $\hat{\sigma}_{i,\text{max}}$ denote the minimum and maximum standard singular values of $Q_i$, respectively, then $\kappa\left(\xTx{Q_i}+\pi I\right)\leq \kappa_\text{max}$ $\forall i=1,\dots,N$,
if $\pi$ is chosen such that
\begin{align}\label{eq:pichoice}
\pi \geq \min_{i\in\lbrace 1,\dots,N\rbrace}\frac{\hat{\sigma}_{i,\text{max}} - \kappa_\text{max}\hat{\sigma}_{i,\text{min}}}{\kappa_\text{max}-1}.
\end{align}

\begin{figure}
\ifArxivFormat
\centering
\includegraphics[scale=1]{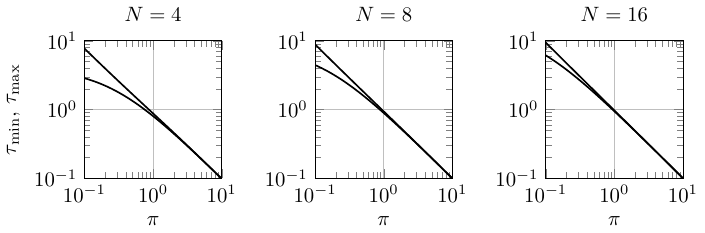}
\else
\input{figure_tau.tex}
\fi
\caption{Minimum and maximum eigenvalues of $T_\pi$ as a function of $\pi$ for different $N$.}\label{fig:tau}
\end{figure}

After obtaining the eigenvectors of $T_\pi$ on line~\ref{alg:HOGSVD:EV}, the indices associated with the isolated HO-CSD subspace, $\mathcal{I}_1$, are determined by
\begin{align}\label{eq:isochoice}
\mathcal{I}_1\eqdef \set{k\in\lbrace 1,\dots,n\rbrace}{\frac{\tau_\text{max}-\tau_k}{\tau_\text{max}-\tau_\text{min}}\leq\epsilon},
\end{align}
where $\tau_k$ is the corresponding eigenvalue of $T_\pi$ and the scalar $\epsilon \geq 0$ is introduced to account for finite machine precision. Note the trade-off between~\eqref{eq:pichoice} and~\eqref{eq:isochoice}: For increasing $\pi$, the difference $\tau_\text{max}-\tau_\text{min}$ rapidly decreases, such as shown in Figure~\ref{fig:tau}. If the difference $\tau_\text{max}-\tau_\text{min}$ is too small, numerical inaccuracies can lead to a wrong selection of directions associated with the isolated HO-CSD subspace. The same problem arises when determining the common HO-CSD subspace.

If $\dim\left(\mathcal{T}_1\lbrace Q_1,\dots,Q_N\rbrace\right)\geq 1$, $T_\pi$ has an eigenvalue that is equal to $\tau_\text{max}$ with geometric multiplicity greater than 1. In this case, the corresponding eigenvectors of $T_\pi$ must be aligned with the shared standard right singular vectors of $Q_1,\dots,Q_N$. Otherwise, it is not guaranteed that the corresponding generalized singular values are either $0$ or $1$. %In~\cite{HOGSVDGITHUB}, this is accomplished by computing the standard SVD of $Q_i Z_{\mathcal{I}_1}$, $i=1,\dots,N$, and identifying the unit singular values analogously to~\eqref{eq:isochoice}.
Section~\ref{sec:computationiso} summarises how this can be accomplished by computing a sequence of standard SVDs.%, and identifying the unit singular values analogously to~\eqref{eq:isochoice}.

Finally, the matrices $\Sigma_i$ and $U_i$ are computed in the loop starting on line~\ref{alg:HOGSVD:loop}. If $\sigma_{i,k}=0$, the algorithm substitutes one of the basis vectors spanning $\ker(\trans{A})$. Note that if $\rank(Q_i)=r_i<n$, then $n-r_i$ generalized singular values are zero and $\dim(\ker(\trans{A}))=n-r_i$. On line~\ref{alg:HOGSVD:V}, the shared matrix of right basis vectors $V$ is obtained without the need for computing the inverse of $R$. 

An upper bound on the algorithm complexity is given by summing the shaded numbers on the right-hand side of Algorithm~\ref{alg:HOGSVD}. The algorithm mainly uses standard routines, such as the QR decomposition or the eigendecomposition, which require roughly $\mathcal{O}\left(M n^2\right)$ floating-point operations. However, forming the matrix $T_\pi$ requires $N$ matrix inversions of size $n\times n$ or $\mathcal{O}\left(N n^3\right)$ floating-point operations. The accumulation of the $N$ inverses can also lead to a non-trivial loss of accuracy. If the full factorization~\eqref{eq:factorization} is not required but only the common or isolated subspace, alternative algorithms exist that compute the common HO-CSD subspace from intersecting the pairwise common HO-GSVD subspaces of $Q_i$ and $Q_{i+1}$ for $i=1,\dots N-1$~\cite{LOANSLIDES,HOGSVDTHESIS}. The pairwise subproblems can be solved by the standard GSVD, which exists as a built-in function in most scientific computing packages, and specialized algorithms exist for large-scale problems~\cite{LARGEGSVD}. According to Theorem~\ref{thm:HOGSVD_GSVD2}, the common subspace algorithm from~\cite{LOANSLIDES} can also be adapted for the isolated HO-CSD subspace.

\subsection{\added{Computing the Isolated Subspace}}\label{sec:computationiso}
It follows from Definition~\ref{def:isolatedsubspace} and Theorem~\ref{thm:HOCSDisolated} that if \[\dim\left(\mathcal{T}_1\lbrace Q_1,\dots,Q_N\rbrace\right)\reqdef n_\text{iso},\] then $T_\pi$ has $n_\text{iso}$ eigenvalues equal to $\tau_\text{max}$ and each of the corresponding eigenvectors can be chosen such that it is a standard right singular vector for each $Q_i$. However, when $n_\text{iso} > 1$ the eigendecomposition of $T_\pi$ will produce an arbitrary set of orthogonal vectors that span $\mathcal{T}_1\lbrace Q_1,\dots,Q_N\rbrace$, but that are not necessarily parallel to the shared right standard singular vectors. By Definition~\ref{def:HOCSD}, the eigenvectors of $T_\pi$ spanning $\mathcal{T}_1\lbrace Q_1,\dots,Q_N\rbrace$ must be aligned with the corresponding standard right singular vectors in order to obtain generalized singular values that are equal to $0$ or $1$.

Given $Z_{\mathcal{I}_1}$ that has been obtained from~\eqref{eq:isochoice} for some $\epsilon>0$, one way to align the columns of $Z_{\mathcal{I}_1}$ is to compute the standard SVDs of $Q_i Z_{\mathcal{I}_1}$ for each $i$, and select those directions associated with standard singular values $\hat{\sigma}_{i,k}$ that satisfy 
\begin{align}\label{eq:isochoice2}
\mathcal{I}_1^i\eqdef \set{k\in\lbrace 1,\dots,n\rbrace}{\hat{\sigma}_{i,k}\geq 1-\tilde{\epsilon}},\qquad i=1,\dots,N,
\end{align}
for some other $\tilde{\epsilon}>0$. However, in the presence of numerical inaccuracies, it is unclear how to choose $\tilde{\epsilon}$ to obtain exactly $n_\text{iso}$ directions from~\eqref{eq:isochoice2}, where $n_\text{iso}$ is determined from~\eqref{eq:isochoice} for a given $\epsilon\geq 0$. Since from~\eqref{eq:assumptionQ},~\eqref{eq:T} and Lemma~\ref{thm:P} it follows that $\twonorm{T_\pi z}$ is maximised if $z$ is parallel to the standard right singular vector of $Q_i$ associated with the largest singular value, it appears natural to order the $Q_i$ by magnitude of $\twonorm{Q_i Z_{\mathcal{I}_1}}$, and then select the standard right singular vector associated with the largest $\twonorm{Q_i Z_{\mathcal{I}_1}}$.

Algorithm~\ref{alg:isosub} computes a sequence of ever-thinner standard SVDs to obtain an aligned basis $W_{\mathcal{I}_1}$ from $Z_{\mathcal{I}_1}$, where $\trans{W}_{\mathcal{I}_1}W_{\mathcal{I}_1}=I$ and the columns of $W_{\mathcal{I}_1}$ span the same subspace as those of $Z_{\mathcal{I}_1}$. In the first iteration, the algorithm selects the class~$i$ that has the maximum amplification in the subspace spanned by the columns of $Z_{\mathcal{I}_1}$, i.e. by comparing $\twonorm{Q_i Z_{\mathcal{I}_1}}$. The corresponding direction $Z_{\mathcal{I}_1}\hat{v}_1$ is assigned to the first column of $W_{\mathcal{I}_1}$. Next, the algorithm selects $n_\text{iso}-1$ remaining directions that are orthogonal to $Z_{\mathcal{I}_1}\hat{v}_1$. Since $Z_{\mathcal{I}_1}$ is orthogonal and at every iteration $\hat{v}_1$ is orthogonal to $\hat{v}_2,\dots,\hat{v}_{n_\text{iso}-k}$, the resulting $W_{\mathcal{I}_1}$ is orthogonal too. Note that the size of $X_k$ decreases at every iteration and that line~\ref{notlast} is not executed at the last iteration. However, given that the first iteration of algorithm Algorithm~\ref{alg:isosub} is of the same worst-case complexity as~\eqref{eq:isochoice}, which amounts to $\mathcal{O}\left(2Mn^2+n^3\right)$ floating point operations. Both methods -- Algorithm~\ref{alg:isosub} as well as~\eqref{eq:isochoice2} -- are implemented in~\cite{HOGSVDGITHUB}, and in Section~\ref{sec:application}, all examples are computed using~\eqref{eq:isochoice2}.
\begin{algorithm}[H]
 \caption{Isolated Subspace Computation}\label{alg:isosub}
 \begin{algorithmic}[1]
    \REQUIRE $Q_1,\dots,Q_N$, $Z_{\mathcal{I}_1}$
    \ENSURE Aligned basis $W_{\mathcal{I}_1}$
    \STATE Initialize $X_0\eqdef Z_{\mathcal{I}_1}$
	\FOR{$k = 0,\dots,n_\text{iso}-1$}
		\STATE Select $p\eqdef\argmax_i \,\twonorm{Q_i X_k}$
		\STATE Obtain the standard right singular vectors $\hat{v}_1,\dots,\hat{v}_{n_\text{iso}-k}$ of $Q_pX_k$
		\STATE Assign $X_k \hat{v}_1$ to column $k+1$ of $W_{\mathcal{I}_1}$
		\STATE Update $X_{k+1}\eqdef X_k\begin{bmatrix}\hat{v}_2 & \dots & \hat{v}_{n_\text{iso}-k}\end{bmatrix}$\label{notlast}
	\ENDFOR
 \end{algorithmic}
\end{algorithm}

\section{\added{Applications}}\label{sec:application}
The standard (HO-)GSVD has already been applied in various fields such as bioinformatics~\cite{HOGSVD,GENES}, medicine~\cite{VACCINE}, acoustics~\cite{HOGSVDSENSORS} or control theory~\cite{MULTIARRAYGSVD}. In practice, the HO-GSVD is used to compare $N$ sets of measurements tabulated in matrices $A_1,\dots,A_N$, where matrix $i$ represents a different organism, class or experiment, for example. The columns of $A_i$ usually represents a sampled coordinate, such as time or position, whereas the rows of $A_i$ class-specific variables that vary along the sampled coordinate.

In the form of the HO-GSVD factorization~\eqref{eq:factorization2}, row $j$ of $A_i$ is represented as a linear combination of the right basis vectors $v_1,\dots,v_n$, which are also weighted by the generalized singular values $\sigma_{i,k}$. In general, the right basis vectors are not orthogonal. However, suppose $A_1,\dots,A_N$ are such that there exists $v\in\R^n$ such that $\xTx{A_i} v\neq 0$ for some $i$ and $\xTx{A_j} v = 0$ for $j\neq i$, i.e. $v$ contributes exclusively to the rows of $A_i$, then, according to Corollary~\ref{thm:HOGSVDisolated}, $v$ will be an eigenvector of $T_\pi$ associated with an eigenvalue equal to $\tau_\text{max}$. Due to the continuity of the eigenvalues of $T_\pi$, it also follows that if $\tau_k\approx\tau_\text{max}$, the corresponding right basis vector is almost exclusively used to represent the rows of $A_i$ (see also~\cite[Ch. 2.3.3]{HOGSVDTHESIS}). Similarly, if there exists $\tilde{v}\in\R^n$ such that $\xTx{A_i} \tilde{v}= \xTx{A_j} \tilde{v}$ for $i,j=1,\dots,N$, then according to Statement~\ref{thm:HOGSVDcommonIV} of Corollary~\ref{thm:HOGSVDcommon}, $\inv{D_{j,\pi}}\tilde{v}$ will be an eigenvector of $T_\pi$ associated with an eigenvalue equal to $\tau_\text{min}$. Among other cases, the condition $\xTx{A_i} \tilde{v}= \xTx{A_j} \tilde{v}$ holds if $A_1,\dots,A_N$ share a singular vector $\tilde{v}$ associated with an identical singular value.

To examine the effect of certain right basis vectors onto the rows of class $i$, $A_i$ can be reconstructed by using a reduced set of right basis vectors, e.g. computing 
\begin{align}\label{eq:Aiiso}
A_{i,\text{iso}}\eqdef\sum_{k\in\mathcal{I}_1} \sigma_{i,k}u_{i,k}\trans{v}_k,
\end{align}
yields the reconstruction of $A_i$ using the right basis vectors that are, in the sense of~\eqref{eq:isochoice}, exclusively used by class $i$. To see the effect of the right basis vectors associated with the common subspace, $A_i$ can be reconstructed by summing over $k\in\mathcal{I}_N$.

\subsection{\added{Numerical Example}}
To illustrate an example application of the HO-GSVD for rank-deficient matrices, consider the CIFAR-10 dataset, which is a collection of images used to evaluate machine learning and computer vision algorithms~\cite{CIFAR10}. The CIFAR-10 dataset provides $6$ batches of $10,000$ $32\times 32$ color images in 10 different classes, and here the rank-deficient HO-GSVD is used to analyse a subset of $N=4$ classes shown in Table~\ref{tab:1}. The following example can be downloaded from~\cite{HOGSVDGITHUB}.
\begin{table}
\caption{Sample matrices extracted from the first batch of the CIFAR-10 dataset. The rows of each $\inR{A_i}{m_i}{n}$ represent vectorised $32\times 32$ pixels large images.}\label{tab:1}
\begin{center}
\begin{tabular}{l c c c c c}
\toprule
 & Class & $m_i$ & $\rank(A_i)$ & $\abs{\mathcal{I}_1^i}$ & $\dim(\mathcal{T}_1\lbrace A_1,\dots,A_4 \rbrace)$\\
\midrule
$A_1$ & Automobile 	& 974 	& 974	& 51\\ 
$A_2$ & Cat 		& 1016 	& 1016	& 92\\ 
$A_3$ & Ship 		& 1025 	& 1025	& 100\\ 
$A_4$ & Truck 		& 981 	& 981	& 57\\[0.5em]
%\hline
$A$ &  				& 3996 	& 3072	& & 300\\
\bottomrule
\end{tabular}
\end{center}
\end{table}

\begin{figure}[]
\centering
\ifArxivFormat
\begin{tabular}{>{\centering}p{.2\textwidth}>{\centering}p{.2\textwidth}>{\centering}p{.2\textwidth}>{\centering\arraybackslash}p{0.25\textwidth}}
\includegraphics[width=.1\textwidth]{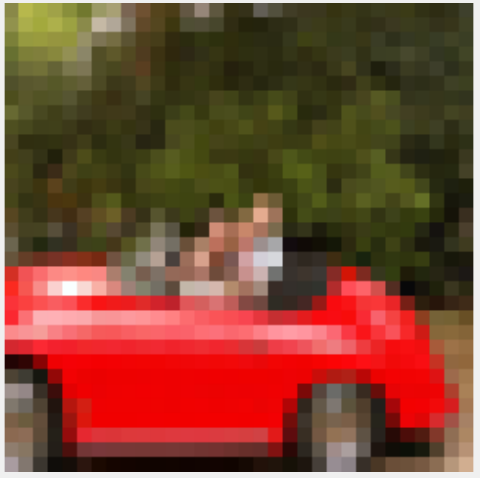} &
\includegraphics[width=.1\textwidth]{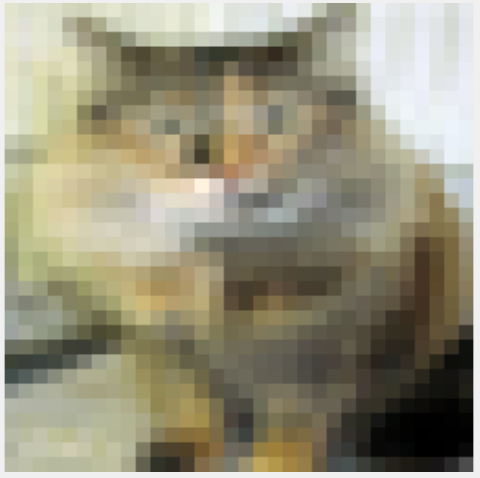} &
\includegraphics[width=.1\textwidth]{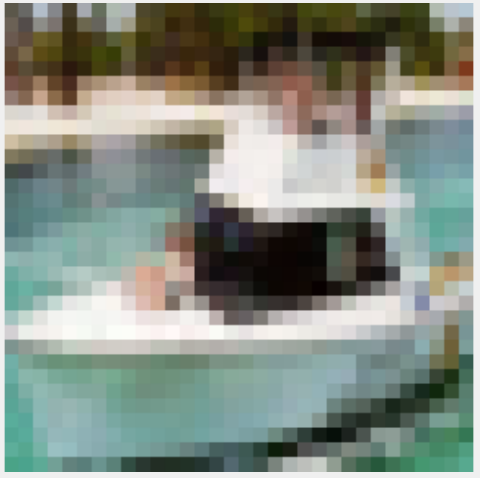} &
\includegraphics[width=.1\textwidth]{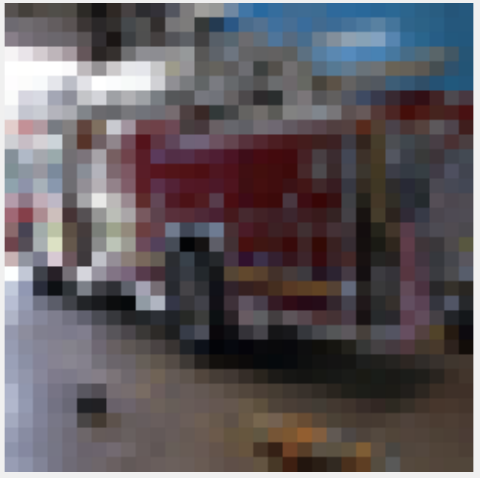} \\[2em]
\includegraphics[width=.1\textwidth]{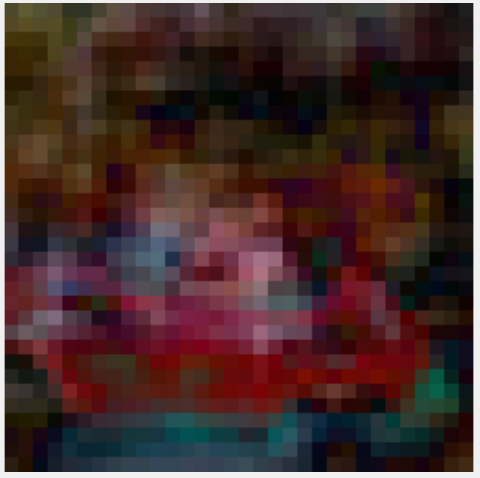} &
\includegraphics[width=.1\textwidth]{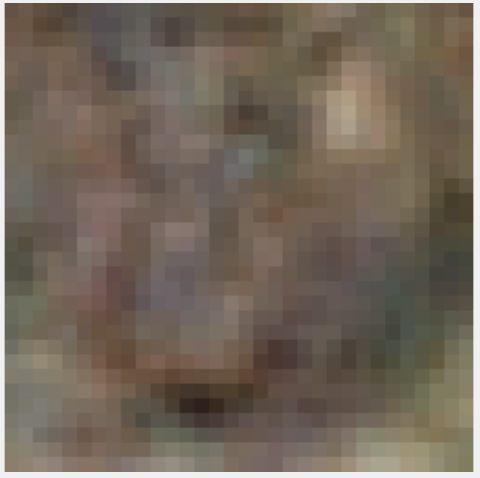} &
\includegraphics[width=.1\textwidth]{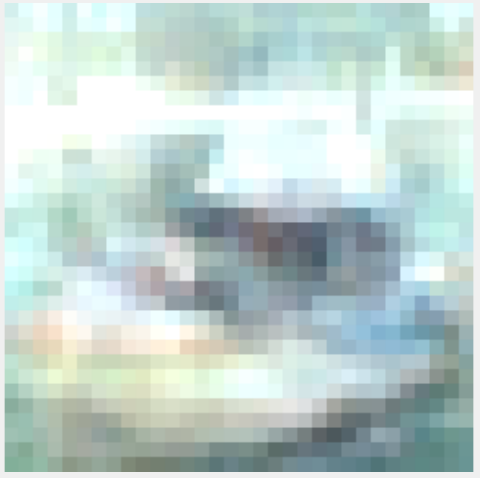} &
\includegraphics[width=.1\textwidth]{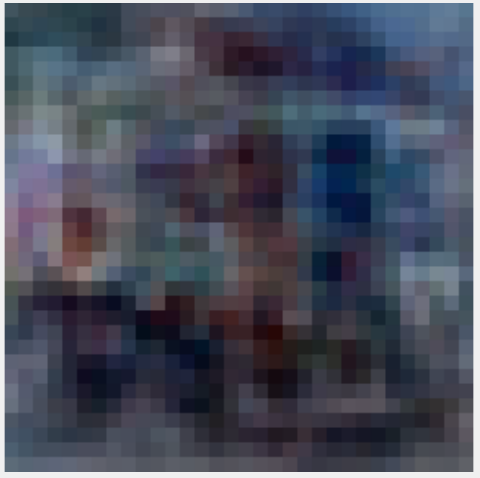} \\[2em]
\includegraphics[width=.1\textwidth]{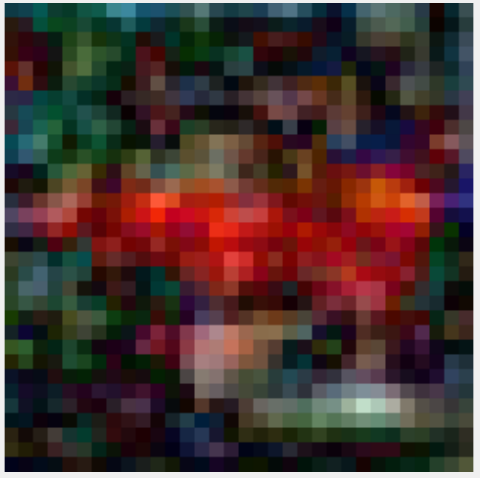} &
\includegraphics[width=.1\textwidth]{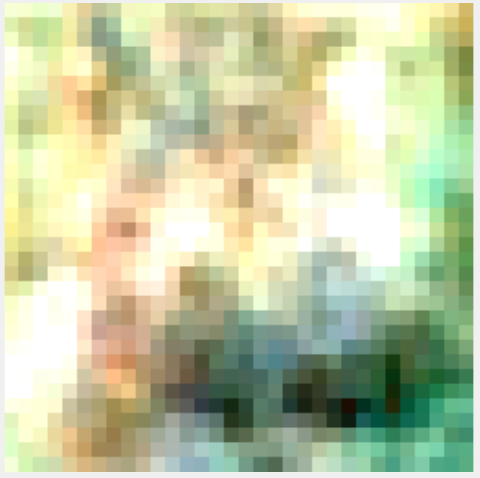} &
\includegraphics[width=.1\textwidth]{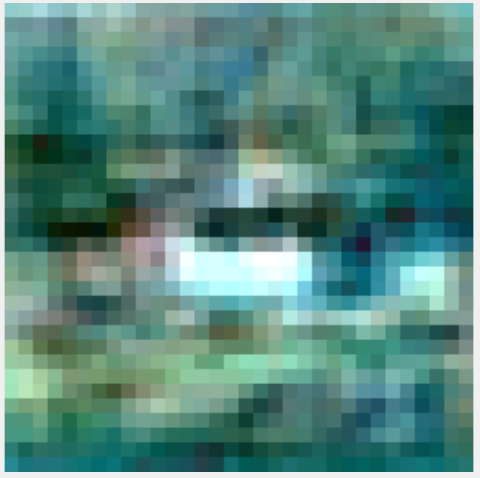} &
\includegraphics[width=.1\textwidth]{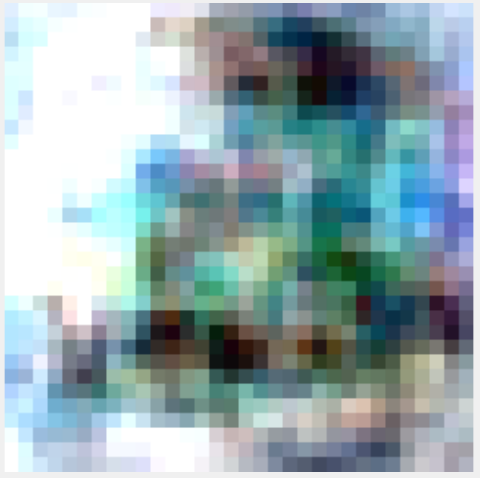}
\end{tabular}
\else
\begin{tabular}{>{\centering}p{.2\textwidth}>{\centering}p{.2\textwidth}>{\centering}p{.2\textwidth}>{\centering\arraybackslash}p{0.25\textwidth}}
\includegraphics[width=.1\textwidth]{figures/Aij_automobile.eps} &
\includegraphics[width=.1\textwidth]{figures/Aij_cat.eps} &
\includegraphics[width=.1\textwidth]{figures/Aij_ship.eps} &
\includegraphics[width=.1\textwidth]{figures/Aij_truck.eps} \\[2em]
\includegraphics[width=.1\textwidth]{figures/Aij_iso_automobile.eps} &
\includegraphics[width=.1\textwidth]{figures/Aij_iso_cat.eps} &
\includegraphics[width=.1\textwidth]{figures/Aij_iso_ship.eps} &
\includegraphics[width=.1\textwidth]{figures/Aij_iso_truck.eps} \\[2em]
\includegraphics[width=.1\textwidth]{figures/v20.eps} &
\includegraphics[width=.1\textwidth]{figures/v82.eps} &
\includegraphics[width=.1\textwidth]{figures/v203.eps} &
\includegraphics[width=.1\textwidth]{figures/v278.eps}^
\end{tabular}
\fi
\caption{First row: Example rows of $A_1,\dots,A_4$ (left to right) reshaped into $32\times 32$ pixels large images. Second row: Moduli of example rows of $A_{1,\text{iso}},\dots,A_{4,\text{iso}}$, where $A_{i,\text{iso}}$ is reconstructed using right basis vectors from the isolated subspace only. Third row: Moduli of isolated right basis vectors that have the largest weight in each image.}\label{fig:images}
\end{figure}

The images are vectorized and grouped in the matrices $\inR{A_i}{m_i}{n}$, where $n=32 \times 32\times 3= 3072$ and $0\leq A_i\leq 1$ (element-wise). Each $A_i$ is such that $r_i\eqdef\rank(A_i)<n$, but the stacked $\inR{A}{M}{n}$ satisfies $M>n$ and $\rank(A)=n$. The first row of Figure~\ref{fig:images} displays row $j_i$ for each class $i$ as an image\footnote{The rows $j_i$ for class $i$ are $j_1=16,j_2=19,j_3=40$ and $j_4=50$, and have been selected to yield an interpretable reconstruction in the isolated subspace.}.

Using the HO-GSVD, the image $j$ of class $i$ can be represented as $\sum_{k} (\trans{e}_ju_{i,k})\sigma_{i,k}\trans{v}_k$, where $e_j$ is a standard basis vector and $\abs{\trans{e}_ju_{i,k}}\leq 1$. The columns $v_k$ of the matrix $\inR{V}{n}{n}$ with $\det(V)\neq 0$ can be interpreted as ``basis images'' for the space of $32 \times 32$ images, and class $i$ uses $r_i$ columns of $V$ to represent its sample images. Note that the columns of $V$ are not orthogonal, and some right basis vectors can therefore ``cancel out'' each other. The third row of Figure~\ref{fig:images} visualises right basis vectors $v_{20}$, $v_{82}$, $v_{203}$ and $v_{278}$, which are all associated with the isolated subspace (see the subsequent paragraphs).
%
%\begin{figure}
%\begin{center}
%\input{figure_U.tex}
%\end{center}
%\caption{Row $1$: The $n=3072$ eigenvalues of $T_\pi$ relative to the bounds $\tau_\text{min}$ and $\tau_\text{min}$. Rows 2--5: Corresponding generalized singular values $\Sigma_i=\diag(\sigma_{i,1},\dots,\sigma_{i,n})$ for classes $i=1,\dots,4$.}\label{fig:Ui}
%\end{figure}

\begin{figure}
\begin{center}
\ifArxivFormat
\centering
\includegraphics[scale=1]{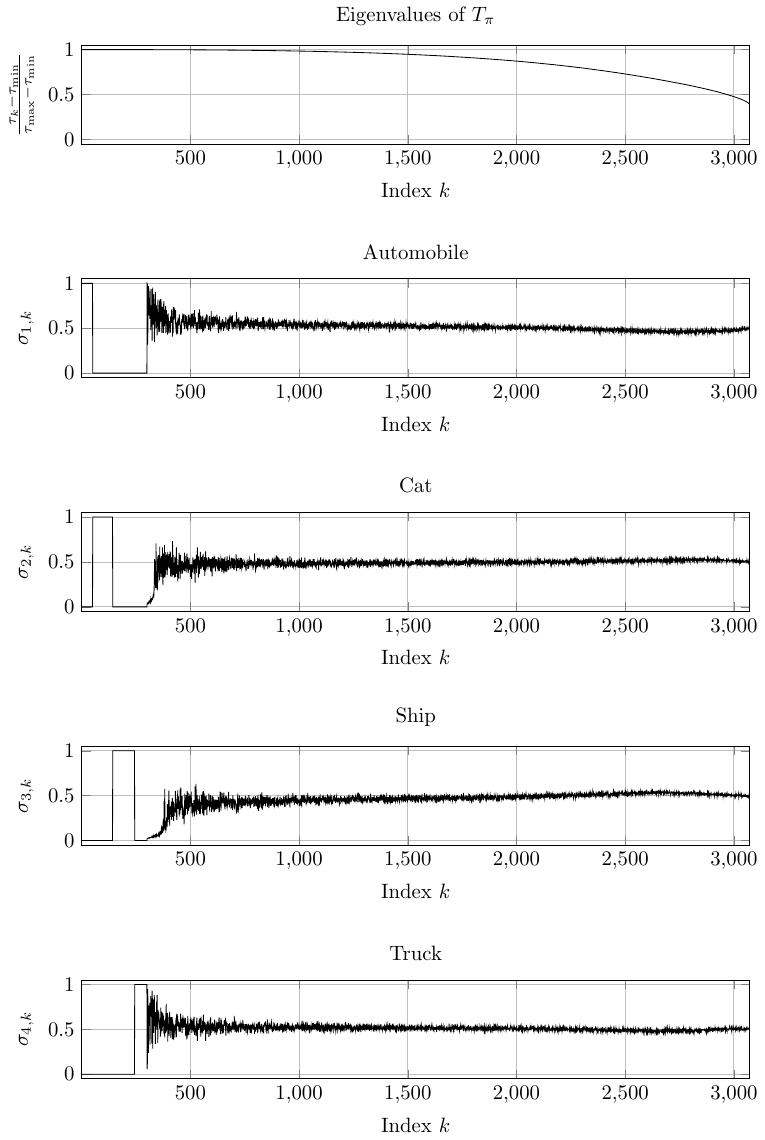}
\else
\input{figure_T.tex}
\fi
\end{center}
\caption{Row $1$: The $n=3072$ eigenvalues of $T_\pi$ relative to the bounds $\tau_\text{min}$ and $\tau_\text{min}$. Rows 2--5: Corresponding generalized singular values $\Sigma_i=\diag(\sigma_{i,1},\dots,\sigma_{i,n})$ for classes $i=1,\dots,4$.}\label{fig:cifar}
\end{figure}

The parameter $\pi$ is chosen as $\pi=1/N=0.25$, which results in $\kappa(\xTx{Q_i}+\pi I)\leq 5$, $\tau_\text{min}=2$ and  $\tau_\text{max}=3.2$. The $n=3072$ eigenvalues $\tau_k$ of $T_\pi$ are shown in the first row of Figure~\ref{fig:cifar}, where $\tau_k$ is displayed relative to $\tau_\text{min}$ and $\tau_\text{max}$ as $(\tau_k-\tau_\text{min})/(\tau_\text{max}-\tau_\text{min})$ sorted in descending order. It can be seen that most eigenvalues are closer to $\tau_\text{max}$ than $\tau_\text{min}$, and that $\tau_k\gg\tau_\text{min}\quad\forall k$, i.e. the common HO-GSVD subspace is empty. Using a tolerance of $\epsilon=10^{-6}$, the dimension of the isolated HO-CSD subspace is estimated as $n_\text{iso}=300$. The number of isolated directions per class is computed from~\eqref{eq:isochoice2} with $\tilde{\epsilon}=\epsilon$, and $\abs{\mathcal{I}_1^i}$ is shown in Table~\ref{tab:1} for each class.

The generalized singular values $\Sigma_i=\diag(\sigma_{i,1},\dots,\sigma_{i,n})$ are shown on the second to fifth row of Figure~\ref{fig:cifar}. For indices $k\in\mathcal{I}_1$ that are associated with the isolated subspaces, the generalized singular values are either 0 or 1. Due to numerical inaccuracies, the separation between $\sigma_{i,k}$, $k\in\mathcal{I}_1$, and $\sigma_{i,j}$, $j\notin\mathcal{I}_1$, is not sharp, i.e. the generalized singular values of the automobile class soar at index $k=301\notin\mathcal{I}_1$ before decreasing at larger indices. Note that even though some $\sigma_{i,k}$ equal $1/\sqrt{N}=0.5$, which is the same magnitude as expected for an index $k$ associated with the common subspace, the common subspace is empty, as can be seen from the first row of Figure~\ref{fig:cifar}.

From Figure~\ref{fig:cifar}, it becomes clear that each class $i$ uses its own subset of isolated basis images as well as $n-n_\text{iso}=2072$ other columns of $V$ to form its $m_i$ samples. Class $i$ can be reconstructed using~\eqref{eq:Aiiso} to obtain $A_{i,\text{iso}}$, which considers indices $k\in\mathcal{I}_1^i$ only. The second row of Figure~\ref{fig:images} shows row $j_i$ of $A_{i,\text{iso}}$, where some degree of resemblance between the original and reconstructed image exists. Examples of right basis vectors are given in the third row of Figure~\ref{fig:images} that shows $v_{20}$, $v_{82}$, $v_{203}$ and $v_{278}$, each of which is associated with the isolated subspace of classes $i=1,\dots,4$. The right basis vectors have been selected by determining those $k$ that maximise $\abs{\trans{e}_{j_i}u_{i,k}}$ for each image $j_i$, i.e. those right basis vectors have a large contribution to image $j_i$. As for the second row of Figure~\ref{fig:images}, it can be seen that the third row of Figure~\ref{fig:images} resembles the original image.

To complement the numerical example, the dataset $A$ is modified in order to artificially introduce a non-empty common subspace. According to Corollary~\ref{thm:HOGSVDcommon}, the common HO-GSVD subspace, $\mathcal{S}_4\lbrace A_1,\dots,A_4\rbrace$, is non-empty iff the condition $\xTx{A_i} \tilde{v}= \xTx{A_j} \tilde{v}$ holds $\forall i,j=1,\dots,4$ and for some $\tilde{v}$, which can be written out as
\begin{align}\label{eq:commoncond}
\left(
\begin{bmatrix}
\trans{a_{i,1}}a_{i,1} & \dots & \trans{a_{i,1}}a_{i,n} \\[3pt]
\vdots & \ddots & \vdots\\[3pt]
\trans{a_{i,n}}a_{i,1} & \dots & \trans{a_{i,n}}a_{i,n}
\end{bmatrix}-
\begin{bmatrix}
\trans{a_{j,1}}a_{j,1} & \dots & \trans{a_{j,1}}a_{j,n} \\[3pt]
\vdots & \ddots & \vdots\\[3pt]
\trans{a_{j,n}}a_{j,1} & \dots & \trans{a_{j,n}}a_{j,n}
\end{bmatrix}\right)\tilde{v}=0,
\end{align}
where $a_{i,k}\in\R^{m_i}$ denotes column $k$ of matrix $A_i$. If $\tilde{v}$ is chosen as $\begin{bmatrix}1 & 0 & \dots & 0\end{bmatrix}^\Tr$, condition~\eqref{eq:commoncond} is tantamount to requiring that $\trans{a_{i,k}}a_{i,1}=\trans{a_{j,k}}a_{j,1}$ $\forall i,j=1,\dots,4$ and for $k=1,\dots,n$, i.e. the projection of column $k$ onto the first column of class $i$ must equal the projection of column $k$ onto the first column of class $j$. Note that condition~\eqref{eq:commoncond} is \emph{not} equivalent to inserting an identical image $x\in\R^n$ in each $A_i$, but a simple way to satisfy~\eqref{eq:commoncond} is to set $a_{i,1}=\begin{bmatrix}1 & 0 & \dots & 0\end{bmatrix}^\Tr$ and zero out the first element of $a_{i,k}$, $k = 1,\dots,n$, for all classes $i=1,\dots,4$. This way the first image of each $A_i$ is replaced with a black square that has one red pixel in the left corner. 

\begin{figure}
\begin{center}
\ifArxivFormat
\centering
\includegraphics[scale=1]{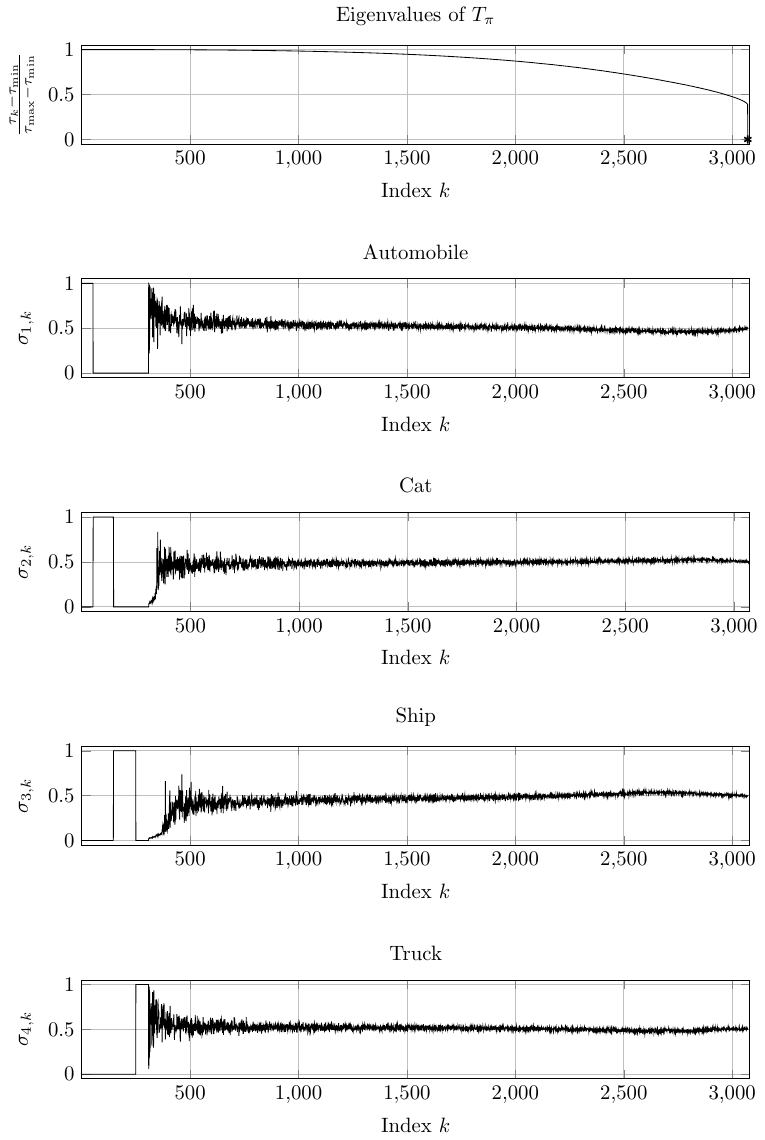}
\else
\input{figure_Tcommon.tex}
\fi
\end{center}
\caption{Row $1$: The $n=3072$ eigenvalues of $T_\pi$ relative to the bounds $\tau_\text{min}$ and $\tau_\text{min}$ for the modified dataset, which has a one-dimensional common subspace associated with index $k=3072$ (marked by an asterisk). Rows 2--5: Corresponding generalized singular values for the modified dataset.}\label{fig:cifaramended}
\end{figure}

The eigenvalues of $T_\pi$ and the generalized singular values for the modified dataset are shown in Figure~\ref{fig:cifaramended}. In the first row of Figure~\ref{fig:cifaramended}, it can be seen that $\tau_k=\tau_\text{min}$ for $k=n$, i.e. the modification successfully introduces a non-empty common subspace. The corresponding generalized singular values equal $1/\sqrt{N}=0.5$ for each class. By construction, the first row (image) of each $A_i$ is orthogonal to all other rows of $A_i$, and therefore aligned with a shared standard right singular vector. The right basis vector associated with the common subspace, $v_n$, is therefore orthogonal to all other basis vectors, which is not the case in general. However, for this example it follows that $v_n$ contributes equally to each of the matrices $A_i$, and for each class $i$, $v_n$ is used to represent the first image only.

\section{Conclusion}
In this paper, we have extended the standard HO-GSVD~\cite{HOGSVD} to accommodate column rank-deficient matrices. By adding the term $\pi\xTx{A}$ to each of the quotient terms $D_{i,\pi}=\xTx{A_i}+\pi\xTx{A}$, we shifted their eigenvalues and bounded them away from zero. This allowed the full-rank requirement on each $A_i$ to be omitted and to extend the HO-GSVD with the notion of isolated subspaces.

The choice of adding a multiple of $\xTx{A}$ was motivated by the relationship between $S_{\pi}$ and $T_{\pi}$, which yielded the same relationship than in~\cite{HOGSVD} for $\pi =0$.  We bounded the eigenvalues of $T_\pi$ and showed that the extremal eigenvalues are attained iff the corresponding eigenvectors are standard right singular vectors for each $Q_i$ associated with a particular singular value. This led to the definition of the common and isolated HO-CSD (HO-GSVD) subspaces. In Appendix~\ref{app:T}, we also showed that if the $Q_i$ share a right singular vector $v$ associated with a zero singular value for $P$ matrices $Q_i$ and with an identical singular value for the other $N-P$ matrices $Q_j$, then $T_\pi$ will have a particular eigenvalue $\tau(P)$ associated with the eigenvector $v$. Future research could investigate whether a biconditional (``iff$\,$") connection holds.

The parameter $\pi$ was assumed to be positive, but otherwise left unspecified. The common and isolated HO-CSD and HO-GSVD subspaces are identified irrespective of the value of $\pi$, but other right basis vectors can be rotated for increasing values of $\pi$, and we have investigated the behavior of these vectors for $\pi\rightarrow 0$ and $\pi\rightarrow\infty$. \added{For $\pi\rightarrow 0$, the eigenvectors of $T_\pi$ are solely determined by the rank-deficient $Q_i$, whereas for $\pi\rightarrow\infty$, the eigenvectors of $T_\pi$ converge to those of the mean of symmetrized products of $\xTx{Q_i}$ and $\xTx{Q_j}$.}

\added{In addition, the choice of $\pi$ also affects the condition number of $\xTx{Q_i}+\pi I$, which must be inverted to obtain $T_\pi$, as well as the range of admissible eigenvalues of $T_\pi$, $\tau_\text{max}-\tau_\text{min}$. A large $\pi$ improves the conditioning of $\xTx{Q_i}+\pi I$, but also tightens the range of eigenvalues, which can lead to a wrong estimate of the common or isolated subspaces in the presence of numerical errors.} The optimal choice of $\pi$ remains unclear and future research could investigate the role of the weight $\pi$.

The majority of our developments were based on the HO-CSD. Using the QR factorization of $A=[\trans{A}_1,\dots,\trans{A}_N\trans{]}$, each $A_i$ was represented as $A_i=Q_i R$ and our findings were developed for $Q_1,\dots,Q_N$, which required $A$ to have full column rank. \added{For rank-deficient $A$, we have shown how $A$ can be padded using an additional matrix $A_{N+1}$ to guarantee that $\det(R)\neq 0$. The properties of $A_1,\dots,A_N$ were inferred from the HO-CSD, which allowed to avoid computing the inverse of a potentially ill-conditioned $R$, but a full factorization still requires to invert the terms $\xTx{Q_i}+\pi I$, which can lead to significant numerical errors for large-scale matrices. Future research could focus on finding a possibly iterative algorithm that finds the eigenvectors of $T_\pi$ without the need for inverting the terms $\xTx{Q_i}+\pi I$.}

For the full-rank case, it has been shown that the common subspace can be found using a variational approach~\cite{LOANSLIDES} and that the vectors $v$ spanning the common subspace are stationary vectors for the function $f_\pi(v)$~\eqref{eq:f} with $\pi=0$. We have shown that the same holds for $\pi>0$. It remains unclear how the right basis vectors, which are not in the common or isolated subspaces, are related to $f_\pi(v)$ and whether an eventual connection would lead to a particular choice of the parameter $\pi$.
\bibliographystyle{siamplain}
\bibliography{master_bib_abbrev}

\appendix
\section{Relation between $S_{\pi}$ and $T_{\pi}$}\label{app:SandT}
Let $D_{i,\pi}=\xTx{A_i}+\pi \xTx{A}$ and define $K_i\eqdef\xTx{Q_i}+\pi I$.
Using~\eqref{eq:Di2}, the matrix $S_{\pi}$ is written as
\ifArxivFormat
\begin{align*}
S_\pi=\frac{1}{N(N-1)}\sum_{i=1}^N\sum_{j=i+1}^N\left(D_{i,\pi}\inv{D_{j,\pi}}+D_{j,\pi}\inv{D_{i,\pi}}\right)=\frac{1}{N(N-1)}\trans{R}\left(\sum_{i=1}^N\sum_{j=i+1}^N\!K_i\inv{K_j}\!\!+K_j\inv{K_i}\right)\tinv{R},
\end{align*}
\else
\begin{align*}
S_\pi&=\frac{1}{N(N-1)}\sum_{i=1}^N\sum_{j=i+1}^N\left(D_{i,\pi}\inv{D_{j,\pi}}+D_{j,\pi}\inv{D_{i,\pi}}\right)\\
&=\frac{1}{N(N-1)}\trans{R}\left(\sum_{i=1}^N\sum_{j=i+1}^N\!K_i\inv{K_j}\!\!+K_j\inv{K_i}\right)\tinv{R},
\end{align*}
\fi
so that by considering $\sum_{i=1}^N K_i=(1+\pi N)\I$
\ifArxivFormat
\begin{align*}
\tinv{R}S_\pi\trans{R} =
\frac{1}{N(N\!-\!1)}\sum_{i=1}^N\sum_{j=i+1}^N\!\!K_i\inv{K_j}\!\!+K_j\inv{K_i}
\!=\!\frac{1}{N(N\!-\!1)}\sum_{i=1}^N\!K_i\sum_{j=1}^N\!\inv{K_j}-\frac{1}{N\!-\!1}\I=\frac{1}{N-1}\left((1+\pi N)T_\pi-\I\right).
\end{align*}
\else
\begin{align*}
\tinv{R}S_\pi\trans{R} &=
\frac{1}{N(N\!-\!1)}\sum_{i=1}^N\sum_{j=i+1}^N\!\!K_i\inv{K_j}\!\!+K_j\inv{K_i}
=\frac{1}{N(N\!-\!1)}\sum_{i=1}^N\!K_i\sum_{j=1}^N\!\inv{K_j}-\frac{1}{N\!-\!1}\I\\&=\frac{1}{N-1}\left((1+\pi N)T_\pi-\I\right).
\end{align*}
\fi
\section{Intermediate Eigenvalues of $T_{\pi}$}\label{app:T}
If there exists a vector $t$ with $\twonorm{t}=1$ in the nullspace of $P$ matrices $Q_j$, but in the range of all other $Q_i$ with index $i\in\mathcal{R}$, then the inequalities~\eqref{eq:ineqTa}-\eqref{eq:ineqTb} can be reformulated as
\begin{subequations}\label{eq:ineqTt}\begin{align}
\trans{t}T_\pi t &= \frac{1}{N}\sum_{i=1}^N\trans{t}\invbr{\xTx{Q_i}+\pi I}t= \frac{P}{\pi N}+\frac{1}{N}\sum_{i\in\mathcal{R}}\trans{t}\invbr{\xTx{Q_i}+\pi I}t\nonumber\\
&\geq \frac{P}{\pi N}+\frac{1}{N}\sum_{i\in\mathcal{R}}\frac{1}{\trans{t}(\xTx{Q_i}+\pi I)t}\label{eq:ineqTta}\\
&\geq \frac{P}{\pi N}+\frac{N-P}{N}\frac{N-P}{\pi (N-P)+\underbrace{\sum_{i\in\mathcal{R}}\trans{t}(\xTx{Q_i})t}_{=1}} = \frac{P(1-\pi N)+\pi N^2}{\pi N(1+\pi (N-P))}.\label{eq:ineqTtb}
\end{align}\end{subequations}

The term on the right-hand side of~\eqref{eq:ineqTtb} corresponds to the minimum and maximum eigenvalues of $T_\pi$ for $P=0$ and $P=N-1$, respectively. If there exists a shared vector $t$ in the nullspace of $P$ matrices $Q_j$, but in the range of all other $Q_i$, then an eigenvalue of $T_\pi$ will be equal to the corresponding value on the right-hand side of~\eqref{eq:ineqTtb}. Note that~\eqref{eq:ineqTt} does \emph{not} prove the converse.

\section{The arithmetic mean of amplification quotients}\label{app:rayleigh}
The HO-GSVD is related to the function $f_\pi(v)$~\eqref{eq:f}, which can be simplified using the stacked QR decomposition~\eqref{eq:A} as
\begin{align*}
g_\pi(z)=\frac{1}{N(N-1)}\sum_{i=1}^{N-1}\sum_{j=1}^N\left(
\frac{\trans{z}(\xTx{Q_i}+\pi I)z}{\trans{z}(\xTx{Q_j}+\pi I)z}+\frac{\trans{z}(\xTx{Q_j}+\pi I)z}{\trans{z}(\xTx{Q_i}+\pi I)z}\right)\geq 1,
\end{align*}
where $z\eqdef R v$. The gradient $\nabla g_\pi(z)$ of $g_\pi(z)$ is given by
\ifArxivFormat
\begin{align*}
\nabla g_\pi(z) \eqdef \frac{1}{N(N-1)}\sum_{i=1}^{N-1}\sum_{j=1}^N\bigg(&
\frac{1}{\trans{z}W_{j,\pi}z}\left(W_{i,\pi} z - \frac{\trans{z}W_{i,\pi}z}{\trans{z}W_{j,\pi}z} W_{j,\pi} z\right)+\frac{1}{\trans{z}W_{i,\pi}z}\left(W_{j,\pi} z - \frac{\trans{z}W_{i,\pi}z}{\trans{z}W_{j,\pi}z} W_{i,\pi} z\right)
\bigg),
\end{align*}
\else
\begin{align*}
\nabla g_\pi(z) = \frac{1}{N(N-1)}\sum_{i=1}^{N-1}\sum_{j=1}^N\bigg(&
\frac{1}{\trans{z}W_{j,\pi}z}\left(W_{i,\pi} z - \frac{\trans{z}W_{i,\pi}z}{\trans{z}W_{j,\pi}z} W_{j,\pi} z\right)\\
&+\frac{1}{\trans{z}W_{i,\pi}z}\left(W_{j,\pi} z - \frac{\trans{z}W_{j,\pi}z}{\trans{z}W_{i,\pi}z} W_{i,\pi} z\right)
\bigg),
\end{align*}
\fi
where $W_{i,\pi}\eqdef\xTx{Q_i}+\pi I$. To show that $\nabla g_\pi(z)=0$ for $z\in\mathcal{T}_{N}\lbrace Q_1,\dots,Q_N\rbrace$ or $z\in\mathcal{T}_{1}\lbrace Q_1,\dots,Q_N\rbrace$, note that $z$ must be a right singular vector for each $Q_i$. It follows that $W_{i,\pi}z = (\sigma_{i,1}+\pi)z$ and
\begin{align*}\label{eq:deltag2}
W_{i,\pi} z - \frac{\trans{z}W_{i,\pi}z}{\trans{z}W_{j,\pi}z} W_{j,\pi} z =
(\sigma_{i,1}+\pi)z - \frac{\sigma_{i,1}+\pi}{\sigma_{j,1}+\pi}(\sigma_{j,1}+\pi)z=0,
\end{align*}
so that $\nabla g_\pi(z)=0$ if $z\in\mathcal{T}_{N}\lbrace Q_1,\dots,Q_N\rbrace$ or $z\in\mathcal{T}_{1}\lbrace Q_1,\dots,Q_N\rbrace$ for any value of $\pi$. The proof is analogous for the HO-GSVD subspaces.
\end{document}